%% file: Revision1.tex
\documentclass[11pt,a4paper]{article}

\usepackage[in]{fullpage}
\usepackage{bm,amsmath,amssymb,amsthm,algorithm,algorithmic,subfigure,graphicx,mathabx,color}
\newtheorem{theorem}{Theorem}[section]
\newtheorem{definition}{Definition}[section]

\def\la{\left\langle}
\def\ra{\right\rangle}
\def\lb{\left(}
\def\rb{\right)}
\def\lcb{\left\{}
\def\rcb{\right\}}

\def\lsb{\left[}
\def\rsb{\right]}

\def\A{\mathcal{A}}
\def\C{\mathbb{C}}

\def\P{\mathcal{P}}
\def\T{\mathcal{T}}
\def\H{\mathcal{H}}

\def\M{\mathcal{M}}
\def\R{\mathbb{R}}
\def\C{\mathbb{C}}

\def\N{\mathcal{N}}
\def\CS{\mathcal{C}}

\def\BZ{\bm{Z}}

\def\BL{\bm{L}}
\def\BU{\bm{U}}
\def\BV{\bm{V}}
\def\BB{\bm{B}}
\def\BC{\bm{C}}
\def\BD{\bm{D}}
\def\BP{\bm{P}}

\def\BSS{\bm{S}}
\def\BM{\bm{M}}
\def\BQ{\bm{Q}}
\def\BX{\bm{X}}
\def\BA{\bm{A}}
\def\BB{\bm{B}}
\def\BW{\bm{W}}
\def\BC{\bm{C}}
\def\BD{\bm{D}}
\def\BG{\bm{G}}
\def\bgam{\bm{\gamma}}
\def\by{\bm{X}}
\def\BR{\bm{R}}
\def\BS{\bm{\Sigma}}
\def\BI{\bm{I}}

\def\ba{\bm{a}}
\def\bx{\bm{x}}
\def\be{\bm{e}}
\def\bz{\bm{z}}
\def\by{\bm{y}}
\def\bb{\bm{b}}

\def\bg{\bm{g}}

\def\bu{\bm{u}}

\def\bzero{\bm{0}}

\DeclareMathOperator{\rank}{rank}

\DeclareMathOperator{\trace}{trace}

\newcommand\numberthis{\addtocounter{equation}{1}\tag{\theequation}}

\title{Exploiting the structure effectively and efficiently in low rank matrix recovery\footnotetext{Authors  are listed alphabetically.}}
\author{Jian-Feng Cai\thanks{Department of Mathematics, Hong Kong University of Science and Technology, Clear Water Bay, Kowloon, Hong Kong SAR, China.}
\and Ke Wei\thanks{School of Data Science, Fudan University, Shanghai, China.}
}

\begin{document}
\maketitle

\begin{abstract}
Low rank model arises from a wide range of applications, including machine learning, signal processing, computer algebra, computer vision, and imaging science. Low rank matrix recovery is about reconstructing a low rank matrix from incomplete measurements. 
In this survey we review  recent developments on low rank matrix recovery, focusing on three typical scenarios: matrix sensing, matrix completion and phase retrieval. An overview of effective and efficient approaches for the problem is given, 
including nuclear norm minimization, projected gradient descent based on matrix factorization, and Riemannian optimization based on the embedded manifold of low rank matrices.
Numerical recipes of different approaches are emphasized while accompanied by  the corresponding theoretical recovery guarantees.\end{abstract}

\input{intro}
\input{convex}
\input{pgd}
\input{manifold}
\input{conclusion}


\bibliographystyle{abbrv}
\bibliography{ref}
\end{document}

%% file: intro.tex
\section{Introduction}
Reconstructing a low rank matrix from incomplete measurements, typically referred to as {\em low rank matrix recovery}, has received extensive investigations during the last decade. 
For conciseness, consider an $n$ by $n$ {\em real and square} matrix $\bm{X}$ which is unknown, and assume $\rank(\BX)=r\ll n$. 
 Let $\A:\mathbb{R}^{n\times n}\mapsto\mathbb{R}^{m}$ be a linear operator from $n\times n$ matrices to $m$-dimensional vectors, which can be defined explicitly as
\begin{equation}\label{eq:linop}
\mathcal{A}(\bm{Z})
=\begin{bmatrix}
\langle \bm{A}_{1},\bm{Z}\rangle\cr
\langle \bm{A}_{2},\bm{Z}\rangle\cr
\vdots\cr
\langle \bm{A}_{m},\bm{Z}\rangle
\end{bmatrix}
\end{equation}
via a set of measurement matrices $\lcb\bm{A}_{\ell}\rcb_{\ell=1}^m\subset\mathbb{R}^{n\times n}$, where $\la\BA_\ell,\BZ\ra=\trace(\BA_\ell^\top\BZ)$ denotes the inner product between $\BA_\ell$ and $\BZ$.
\begin{quote}{\em The goal in low rank matrix recovery is to reconstruct  $\BX$ from $m\ll n^2$ linear measurements of the form $\by=\A(\BX)$}.
\end{quote} This is an ill-posed problem without assuming any structure on $\BX$ since there are more unknowns than equations. However, noticing that the number of degrees of freedom in an $n$ by $n$ rank $r$ matrix is $(2n-r)r$  \cite{Van:SIOPT:13} which can be much smaller than $n^2$ provided $r$ is small, it is reasonable to expect to reconstruct a low rank matrix from fewer than $n^2$ measurements. Moreover, many  effective and effcient approaches have been developed to target low rank matrix recovery, which will be our focus in this review article.

Low rank matrix recovery arises frequently in many research areas of  science and engineering, for example, machine learning, signal processing, computer algebra, computer vision, imaging science, control, and bioinformatics; see \cite{TK:IJCV:92,KV:SPM:96,KP:BI:07,ABB:PNAS:00,JHSX:SIIMS:11,LHKR:TOG:10,FHB:ACC:03,ARR:TIT:14,CSV:CPAM:13,LS:IP:15} and references therein. In these applications, the target of interest is either low rank itself or exhibits a low rank structure after some linear or nonlinear transformations. Also, it is often the case that different applications correspond to different sorts of measurement matrices. In this survey, we will restrict our attention mostly to the following three different scenarios.
\paragraph{Matrix sensing} In this situation,  each measurement matrix $\BA_\ell$ is usually a dense matrix without a particular simple structure, for example $\BA_\ell$ has i.i.d random Gaussian entries. An important application scenario is quantum tomography where one tries to reconstruct an unknown quantum state from experimental data  \cite{kliesch2017guaranteed}.
The state of a quantum system in quantum mechanics can often be described by a low rank matrix while the measurement matrices  are tensor products of Pauli
matrices \cite{Gro:TIT:11}.
\paragraph{Matrix completion} 
The problem here is essentially about completing a low rank matrix from partial observed entries of the matrix. Thus, each matrix measurement matrix $\BA_\ell$ has the form $\BA_\ell = \be_{i_\ell}\be_{j_\ell}^\top$, where $\be_k$ ($k=i_\ell,~j_\ell$) denotes the vector with only one nonzero entry equal to $1$ in the $k$-th coordinate. Let  $\Omega=\{(i_\ell,j_\ell),~\ell=1,\cdots,m\}$ be a set of indices corresponding to the observed entries of an unknown matrix.  The linear operator $\A$ is usually replaced by $\P_{\Omega}$ in matrix completion, where $\P_{\Omega}$ is the associated sampling operator which acquires only the entries indexed by $\Omega$. A well-known application  of matrix completion is in recommendation system \cite{GNOT1992}, where the task is to infer missing ratings given observed ones.  Since a user's preference is typically determined by a few factors, the rating matrix in a recommendation system is approximately low rank.

\paragraph{Phase retrieval}
In phase retrieval, one would like to reconstruct an object from a set of  magnitude or phaseless measurements.  More precisely, letting $\bx\in\R^n$ be an unknown vector, the task in phase retrieval is to reconstruct it from the phaseless measurements $\by$ given by 
\begin{align*}\by=|\BA\bx|^2,\numberthis\label{eq:A}
\end{align*} where $\BA$ is an $m\times n$ matrix. 
 Phase retrieval has found many important applications in imaging problems such as  X-ray crystallography, electron microscopy, diffractive imaging, and astronomical imaging \cite{Ha93a,Buetal07,Mietal08}. Moreover, it can be cast as a low rank matrix recovery problem. To see this, define the {\em rank one} matrix $\BX=\bx\bx^\top$ and let $\ba_\ell^\top$ be the $\ell$-th row of $\BA$. Then, a simple algebra yields that 
\begin{align*}
y_\ell = |\ba^\top_\ell\bx|^2 = \la\ba_\ell\ba_\ell^\top,\BX\ra.
\end{align*}
Noticing the one to one correspondence between $\bx$ and $\BX$, one can easily see that phase retrieval is indeed a rank one matrix recovery problem, where each measurement matrix is given by $\BA_\ell=\ba_\ell\ba_\ell^\top$.\\

From the pioneering work in \cite{CR:FoCM:09,RFP:SIREV:10}, significant progress has been made on low rank matrix recovery. 
 In this article, we would like to outline some basic ideas behind various effective and efficient approaches for low rank matrix recovery, especially on different ways to exploit low rank structures when designing fast algorithms. 
{Additionally, theoretical recovery guarantees for these approaches  will be presented, concerning a question of central importance in low rank matrix recovery:
\begin{quote}
{\em How many measurements are sufficient for a program to be able to successfully reconstruct a low rank matrix?}
\end{quote}}\noindent
Since there is a large body of literature on this topic, it would be difficult to give an exhaustive survey due to the page limit. Interested readers are recommended to consult the other two review articles \cite{DaRo2016,ChenChi18} for more materials.  
 
\subsection{Notation and organization } 
{Following the notation above, }we use bold face upper letters (e.g. $\bm{Z}$) and bold face lower letters (e.g., $\bm{z}$) to denote matrices and vectors respectively, and use the corresponding normal font letters with subindices for their entries (e.g., $Z_{ij}$ and $z_i$ for entries of $\bm{Z}$ and $\bm{z}$ respectively). In particular, we fix $\BX$ to be the underlying rank-$r$ matrix to be recovered and use $\kappa$ to denote the condition number of $\BX$ defined by $\kappa=\sigma_1(\BX)/\sigma_r(\BX)$. Operators are denoted by calligraphic letters (e.g., $\A$ represents the measurement operator). For a given vector, $\|\cdot\|_p$, $p=1,2,\infty$, denotes its $p$-norm. For a given matrix, $\|\cdot\|_2$ stands for the operator norm, $\|\cdot\|_{2,\infty}$ stands for the maximum of $2$-norms of all rows, $\|\cdot\|_F$ stands for the Frobenius norm, and $\|\cdot\|_{\infty}$ stands for the maximum magnitude of all entries.

The rest of this paper is organized as follows. In Section~\ref{sec:convex} we discuss the theory and algorithms for nuclear norm minimization, which is a convex approach for low rank matrix recovery. In Section~\ref{sec:pgd}, the projected gradient descent algorithm based on matrix factorization is presented with recovery guarantees. Section~\ref{sec:manifold} discusses the approaches based on the embedded manifold of low rank matrices, as well as the extensions to more general low rank matrix recovery
problems. We conclude this survey in Section~\ref{sec:conclusion} with a brief discussion.

%% file: convex.tex
\section{Convex approach: Nuclear norm minimization }\label{sec:convex}

Since we are interested in recovering a low rank  matrix $\BX$ from an underdetermined linear system $\by=\A(\BX)$, it is natural to seek the lowest rank matrix consistent with the measurements, which can be formally expressed as 
\begin{equation}\label{eq:rankmin}
\min \rank(\bm{Z})
\quad\mbox{subject to}\quad 
\mathcal{A}(\bm{Z}) = \bm{y}.
\end{equation}
Evidently as long as $\mathcal{A}$ is injective on the set of  matrices of rank at most $r$, $\BX$ will be the unique solution to \eqref{eq:rankmin}. Indeed, 
it has been shown that if $\A$ consists of $m\ge 4nr-4r^2$ generic measurement matrices, then $\A$ will be injective; see \cite{CRWX:AMSA:18,zhiqiang17}  for more details. Despite this, the rank minimization problem is known to be NP-hard and computationally intractable  since it is an extension of the $\ell_0$-minimization problem in compressed sensing \cite{Don:TIT:06,CRT:TIT:06}.

One of the mostly studied approaches in low rank matrix recovery is to replace the rank of $\BZ$ with its nuclear norm $\|\BZ\|_*$ and then solve the following convex relaxation problem:
\begin{equation}\label{eq:nucmin}
\min \|\bm{Z}\|_*
\quad\mbox{subject to}\quad 
\mathcal{A}(\bm{Z}) = \bm{y},
\end{equation}
where the nuclear norm of $\BZ$ is defined as the sum of its singular values, 
$
\|\BZ\|_*=\sum_{i=1}^n\sigma_i(\BZ).
$
It can be shown that the unit nuclear norm ball $\{\BZ~|~\|\BZ\|_*\leq 1\}$ is the convex hull of rank one matrices with unit Frobenius norm \cite{RFP:SIREV:10}. Therefore, nuclear norm minimization is well aligned with $\ell_1$-minimization for compressed sensing where the unit $\ell_1$ ball is the convex hull of one sparse vectors with uint $\ell_2$-norm. Moreover, nuclear norm minimization can be further cast as a semidefinite programming \cite{RFP:SIREV:10} and we can use the off-the-shelf software packages to solve it \cite{TTT:OMS:99}.

\subsection{Recovery guarantees of nuclear norm minimization}
Since in an $n\times n$ rank $r$ matrix the number of degrees of freedom is $(2n-r)r$, the information-theoretic minimum for the necessary  number of measurements $m$ is $O(nr)$.
In this subsection, we investigate the sufficient number of measurements for nuclear norm minimization to achieve a successful recovery of  the underlying low rank matrix $\BX$ under the three measurement models mentioned in the introduction. 
\subsubsection{Matrix sensing}
The guarantee analysis for matrix sensing is typically based on the notion of restricted isometry property, which was originally developed for compressed sensing in \cite{CT:TIT:05} and was extended to low rank matrix recovery in \cite{RFP:SIREV:10}.
\begin{definition}[Restricted Isometry Property (RIP)]
Let $\mathcal{A}$ be a linear operator from $n\times n$ matrices to vectors of length $m$. For any integer $0<r<n$, we say $\mathcal{A}$ satisfies the restricted isometry property if there exists a constant $\delta_r\in(0,1)$ such that
\begin{equation}\label{eq:RIP}
(1-\delta_r)\|\bm{Z}\|_F^2\leq\|\mathcal{A}(\bm{Z})\|_2^2\leq(1+\delta_r)\|\bm{Z}\|_F^2
\end{equation}
holds for any matrix $\BZ$ of rank at most $r$.
\end{definition}
If each measurement matrix $\BA_\ell$ has i.i.d Gaussian entries of mean $0$ and variance $1/m$, then with high probability $\A$ satisfies the RIP with a small constant provided\footnote{The notation $m\gtrsim nr\log n$ means there exists an absolute constant $C>0$ such that $m\geq Cnr\log n$.} $m\gtrsim nr\log n$ \cite{RFP:SIREV:10}. This sampling complexity was subsequently sharpened to $m\gtrsim nr$ in \cite{CP:XXX:10}. 
For quantum tomography where each involved measurement matrix is a tensor product of Pauli matrices, the RIP was established in \cite{YKL11} for $m\geq nr\log^6n.$
When $\delta_{2r}<1$, it is easy to see that $\A$ is an injective operator on matrices of rank at most $r$ and hence $\BX$ is the unique rank $r$ solution to the rank minimization problem. Moreover, the theoretical recovery guarantee of nuclear norm minimization can be established in terms of the RIP.
\begin{theorem}[\cite{RFP:SIREV:10}]
Assume $\mathcal{A}$ satisfies the RIP with constant $\delta_{5r}<c$ for some small universal constant $c>0$. Then the underlying rank $r$ matrix $\bm{X}$ is the unique solution to \eqref{eq:nucmin}. 
\end{theorem}

\subsubsection{Matrix completion}
The RIP states that the sensing operator is approximately isometry when being restricted to low rank matrices. However, this is not true for matrix completion. 
Recall that $\Omega$ is a subset of indices corresponding to the observed entries and $\P_\Omega$ (alias of $\A$ in matrix completion) denotes the associated sampling operator. We can construct a rank-$1$ matrix $\BZ$ with only one nonzero entry (e.g., equal to $1$) outside of $\Omega$. Then it is trivial to see that $\P_\Omega(\BZ)=\bzero$ and the lower bound in \eqref{eq:RIP} will be violated.
Despite this, the recovery guarantee of nuclear norm minimization for matrix completion can be established based on the notion of incoherence. 
\begin{definition}[Incoherence \cite{CR:FoCM:09}]\label{asump:MC1}
Let $\bm{X}\in\mathbb{R}^{n\times n}$ be a rank $r$ matrix with the compact singular value decomposition (SVD) $\bm{X}=\bm{U}\bm{\Sigma}\bm{V}^T$. We say $\bm{X}$ is $\mu_0$-incoherent if there exists a numerical constant $\mu_0>0$ such that
$$
\|\bm{U}\|_{2,\infty}^2\leq\frac{\mu_0 r}{n}
\quad\mbox{and}\quad
\|\bm{V}\|_{2,\infty}^2\leq\frac{\mu_0 r}{n}.
$$
\end{definition}

\begin{theorem}[\cite{CR:FoCM:09,CT:TIT:10,Gro:TIT:11,Rec:JMLR:11,YDChen15}]
Assume $\BX$ is $\mu_0$-incoherent and  each pair of indices $(i_{\ell},j_{\ell})$ in $\Omega$ is sampled independently and uniformly from $\{1,\ldots,n\}\times \{1,\ldots,n\}$ with replacement. Then with high probability $\BX$ is the  unique solution to \eqref{eq:nucmin} provided
\begin{align*}
m\gtrsim nr\log^2 n.
\end{align*}
\end{theorem}

The proof of the above theorem is based on the construction of a dual certification to certify the optimality of the underlying matrix $\BX$. It is worth noting that the assumption that $\BX$ is $\mu_0$-incoherent is closely related to the uniform sampling scheme.  If some important  sampling scheme  is adopted, the incoherence requirement may be removed; see \cite{CBSW15} and references therein.

\subsubsection{Phase retrieval}
In phase retrieval the ground truth solution $\BX=\bx\bx^\top$  is not only low rank but also positive semidefinite. Thus it is reasonable to add one more constraint to \eqref{eq:nucmin} and solve the following trace minimization problem:
\begin{equation}\label{eq:nucminPR}
\min_{\bm{Z}} \mathrm{trace}(\bm{Z})
\quad\mbox{subject to}\quad 
\mathcal{A}(\bm{Z}) = \bm{y},\quad\bm{Z}\succeq\bm{0},
\end{equation}
where  $\|\BZ\|_*$ can be replaced by $\trace(\BZ)$ since they are equal to each other for the class of  positive semidefinite matrices. The above trace minimization program for phase retrieval is widely known as PhaseLift \cite{CSV:CPAM:13}.

To establish the recovery guarantee of  PhaseLift for phase retrieval, we assume each measurement vector $\ba_\ell$ in $\BA$ (see \eqref{eq:A}) is a standard Gaussian vector; that is $\ba_\ell\sim\N(0,\BI_n)$. Unfortunately, the RIP cannot hold for the corresponding linear operator $\A$ here unless $m$ is on the same order as $n^2$, see \cite{CSV:CPAM:13}. That being said, optimal sampling complexity can still be achieved for PhaseLift via the construction of a dual certificate directly based on the Gaussian random sampling model, leading to the following theorem.

\begin{theorem}[\cite{CSV:CPAM:13,CL:FCM:14}]
Assume $\ba_{\ell}\sim\N(0,\BI_n)$ and $\by=|\BA\bx|^2$. Then  with high probability $\BX=\bx\bx^\top$ is the  unique solution to \eqref{eq:nucminPR} provided 
\begin{align*}
m\gtrsim n.
\end{align*}
\end{theorem}

\paragraph{Remark} We have discussed nuclear norm minimization for low rank matrix recovery, but other  
convex optimization methods  are also available \cite{FHB:ACC:03,SS:LT:05,ZCCZ:IPI:12}. Under  the Gaussian measurement model for matrix sensing, more quantitative phase transitions 
for nuclear norm minimization can be characterized based on convex geometry and statistical dimension \cite{CRPW:FoCM:12,ALMT2014}.  
\subsection{Algorithms for nuclear norm minimization}



As stated previously, nuclear norm minimization can be reformulated as a semidefinite programming (SDP)  \cite{CR:FoCM:09,RFP:SIREV:10} which can be further solved  by interior-point methods in polynomial time. 
However, finding the solution by the interior-point methods needs to solve systems of linear equations to compute the Newton direction in each iteration, which can be prohibitive for large $n$ and hence limit the applicability of nuclear norm minimization if an exact solution to \eqref{eq:nucmin} is sought. 

To avoid the huge linear system when computing the  Newton direction in the interior-point methods, many first order algorithms have been developed for certain variants of \eqref{eq:nucmin}. The most challenging part in the design of efficient algorithms is the non-smoothness of the nuclear norm function. Since the nuclear norm function is non-differentiable, its gradient does not exist and one has to use the subgradient, which can be computed as follows:
$$
\partial\|\bm{Z}\|_*=\left\{\bm{U}\bm{V}^T+\bm{W}~|~\bm{Z}=\bm{U}\bm{\Sigma}\bm{V}^T\mbox{ is the compact SVD}, \bm{W}^T\bm{Z}=\bm{Z}\bm{W}=\bm{0},~\|\bm{W}\|_2\leq 1\right\};
$$ 
see  \cite{CR:FoCM:09,CCS:SIOPT:10}.
For a non-smooth convex function, a simple explicit forward subgradient algorithm is not guaranteed to converge until the stepsize is very small. To allow a larger stepsize, implicit backward gradient descent algorithms may be applied. More precisely, to minimize a non-smooth convex function $f(x)$, an implicit backward subgradient descent updates the variables by $x_{k+1}\in x_{k}-\alpha_k\partial f(x_{k+1})$. In order to get $x_{k+1}$ from $x_k$, we need to solve the inclusion equation, whose solution is given by $$x_{k+1}=\arg\min_x\frac12\|x_k-x\|_2^2+\alpha_k f(x).$$ The mapping from $x_k$ to $x_{k+1}$ is known as the proximity operator of $f$ in convex analysis, which plays an important role in many first-order convex optimization algorithms.  Restricting to the nuclear norm function, it turns out that the proximity operator is the well-known singular value thresholding operator \cite{CCS:SIOPT:10}. 
\begin{theorem}[Singular Value Thresholding (SVT)]\label{thm:SVT} Let $\bm{Z}=\sum_{i=1}^{n}\sigma_i\bm{u}_i\bm{v}_i^T$ be the SVD of $\bm{Z}$. Define the singular value thresholding  on $\bm{Z}$ as follows
\begin{align*}
\mathcal{D}_{\tau}(\bm{Z})=\sum_{i=1}^{n}\max\{\sigma_i-\tau,0\}\bm{u}_i\bm{v}_i^T.\numberthis\label{eq:svt_operator}
\end{align*}
Then, $\mathcal{D}_{\tau}$ is the proximity operator of the nuclear norm function, namely,
$$
\mathcal{D}_{\tau}(\bm{Z})=\arg\min_{\bm{Y}}\frac12\|\bm{Z}-\bm{Y}\|_F^2+\tau\|\bm{Z}\|_*.
$$
\end{theorem}

The SVT  operator is not only used  in the backward gradient descent methods but also used in many first order dual algorithms or primal-dual algorithms targeting the variants of \eqref{eq:nucmin}.
 In the following, we give a few examples of such algorithms without providing detailed convergence analysis.

\paragraph{SVT algorithm}
We may approximate nuclear norm minimization by the following one with a strongly convex objective:
\begin{equation}\label{eq:nucfromin}
\min_{\bm{Z}}\|\bm{Z}\|_*+\frac{1}{2\lambda}\|\bm{Z}\|_F^2\quad \mbox{subject to}\quad\mathcal{A}(\bm{Z})=\bm{y}.
\end{equation} 
First note that this approximation is quite accurate. Indeed, it is shown in \cite{ZCCZ:IPI:12} that, with a sufficiently large finite number $\lambda$,  \eqref{eq:nucfromin} has the same solution as \eqref{eq:nucmin}. The  superiority of using \eqref{eq:nucfromin} is that the Lagrangian dual problem of a strongly convex minimization problem is continuously differentiable. Therefore, a gradient ascent algorithm can be applied to the dual problem of \eqref{eq:nucfromin}, known as Uzawa's algorithm. This leads to the following SVT algorithm \cite{CCS:SIOPT:10}:
\begin{equation}\label{eq:SVTalg}
\begin{cases}
\bm{Y}_{k+1}=\bm{Y}_k-\alpha_k\mathcal{A}^*(\mathcal{A}(\bm{Z}_{k})-\bm{y})\cr
\bm{Z}_{k+1}=\mathcal{D}_{\lambda}(\bm{Y}_{k+1}),
\end{cases}
\end{equation}
where $\alpha_k$ is the stepsize. When the stepsize obeys $0<\inf_k\alpha_k\leq\sup_k\alpha_k<2/\|\mathcal{A}\|^2$, it is proved in \cite{CCS:SIOPT:10} that the sequence $\{\bm{Z}_k\}_{k\in\mathbb{N}}$ generated by \eqref{eq:SVTalg} converges to the unique solution of \eqref{eq:nucfromin}. For matrix completion where $\bm{Y}_k$ will be sparse and $\bm{Z}_k$ will be low-rank, the SVT algorithm is capable of  solving large size problems.

\paragraph{Forward-backward splitting}
When there is noise present in the measurements, it is natural to solve a regularization variant of \eqref{eq:nucmin}, 
\begin{equation}\label{eq:nucminnoise}
\min_{\bm{Z}}\frac12\|\mathcal{A}(\bm{Z})-\bm{y}\|_2^2+\lambda\|\bm{Z}\|_*,
\end{equation}
where $\lambda>0$ is a parameter associated with the noise level. 
Since the first term  in the objective function is smooth, a forward explicit gradient descent is good enough to decrease its value. Noting the second term is non-smooth, an implicit backward gradient descent is suitable, which leads to the SVT. Altogether, we obtain the following iteration:
\begin{equation}\label{eq:FBS}
\bm{Z}_{k+1}=\mathcal{D}_{\alpha_k\lambda}(\bm{Z}-\alpha_k\mathcal{A}^*(\mathcal{A}(\bm{Z}_k)-\bm{y})),
\end{equation}
where $\alpha_k>0$ is the stepsize. In \eqref{eq:FBS}, $\bm{Z}-\alpha_k\mathcal{A}^*(\mathcal{A}(\bm{Z}_k)-\bm{y})$ is the forward explicit gradient descent for the first term in \eqref{eq:nucminnoise}, while $\mathcal{D}_{\alpha_k\lambda}$ is the  implicit backward  gradient descent for the second term as shown in Theorem \ref{thm:SVT}. If the stepsize satisfies $0<\inf_k\alpha_k\leq\sup_k\alpha_k<2/\|\mathcal{A}\|^2$, then the sequence $\{\bm{Z}_k\}_{k\in\mathbb{N}}$ generated by \eqref{eq:FBS} converges to a solution of \eqref{eq:nucminnoise}. The forward-backward splitting framework was surveyed in \cite{CW:MMS:05} for general signal processing problems and was studied for low-rank matrix recovery in \cite{MGC:MP:11}.

\paragraph{Alternating direction method of multipliers (ADMM)}
    The alternating direction method of multipliers (ADMM) is an algorithm that attempts to solve a convex optimization problem by breaking it into smaller pieces, each of which will be easier to handle. A key step in ADMM is the splitting of variables, and different splitting schemes lead to different algorithms. We present an example of ADMM for low rank matrix recovery here. By introducing an auxiliary variable, \eqref{eq:nucminnoise} can be rewritten as the following equivalent convex optimization problem: 
\begin{equation}\label{eq:nucminnoise_admm}
\min_{\bm{Z},\bm{Y}}\frac12\|\mathcal{A}(\bm{Z})-\bm{y}\|_2^2+\lambda\|\bm{Y}\|_*+\frac{\mu}{2}\|\bm{Y}-\bm{Z}\|_F^2\quad\mbox{subject to}\quad \bm{Y}=\bm{Z}.
\end{equation}
The associated augmented Lagrangian function is given by
$$
L_{\mu}(\bm{Y},\bm{Z},\bm{\Lambda})=\frac12\|\mathcal{A}(\bm{Z})-\bm{y}\|_2^2+\lambda\|\bm{Y}\|_*+\frac{\mu}{2}\|\bm{Y}-\bm{Z}\|_F^2+\langle \bm{\Lambda},\bm{Y}-\bm{Z} \rangle,
$$
where $\mu>0$ is a parameter and $\bm{\Lambda}$ is the Lagrange multiplier. Then, the application of an augmented Lagrangian method gives
\begin{equation}\label{eq:ALM}
\begin{cases}
(\bm{Y}_{k+1},\bm{Z}_{k+1})=\arg\min_{\bm{Y},\bm{Z}}L_{\mu}(\bm{Y},\bm{Z},\bm{\Lambda}_k)\cr
\bm{\Lambda}_{k+1}=\bm{\Lambda}_k+\alpha_k\mu(\bm{Y}_{k+1}-\bm{Z}_{k+1}),
\end{cases}
\end{equation}
where $\alpha_k$ is the stepsize. Typically, there does not exist a closed solution for the first minimization problem of \eqref{eq:ALM}. A simple yet effective approximation is to use one step of alternating minimization between $\bm{Y}$ and $\bm{Z}$. After simplifying the expressions and applying Theorem \ref{thm:SVT}, we can obtain the following ADMM algorithm
\begin{equation}\label{eq:ADMM}
\begin{cases}
\bm{Z}_{k+1}=\left(\mathcal{A}^*\mathcal{A}+\mu\mathcal{I}\right)^{-1}(\mathcal{A}^*\bm{y}+\mu\bm{Y}_k+\bm{\Lambda}_k)\cr
\bm{Y}_{k+1}=\mathcal{D}_{\lambda/\mu}(\bm{Z}_{k+1}-\bm{\Lambda}_k/\mu)\cr
\bm{\Lambda}_{k+1}=\bm{\Lambda}_k+\alpha_k\mu(\bm{Y}_{k+1}-\bm{Z}_{k+1}).
\end{cases}
\end{equation}
In the matrix completion case, the first step of \eqref{eq:ADMM} has a closed form solution. For other cases, an efficient linear equation solver can be applied. When $\alpha_k\in(0,(\sqrt{5}+1)/2)$, the algorithm is convergent. Several other different ADMM algorithms have been developed for nuclear norm minimization via different splitting schemes; see for example \cite{BPC:FTML:11,CHY:IMANA:12,TY:SIOPT:11,LCWM:Arxiv:10}. \\

The SVT operator is a key ingredient for many other   algorithms, see \cite{LST:MP:12}  for  an implementable proximal point algorithmic
framework on nuclear norm minimization. Actually, there exists a vast literature on soft-thresholding based algorithms for $\ell_1$-norm minimization in compressed sensing, and these algorithms can be easily adapted  to low-rank matrix recovery after we replace the vector soft-thresholding  operator by the SVT operator. For example, we can adapt FISTA \cite{BT:SIIMS:09} for $\ell_1$-minimization to accelerate the forward-backward splitting algorithm mentioned above \cite{TY:PJO:10}.  In most of the SVT-based algorithms, the main computational cost lies in  the evaluation of  $\mathcal{D}_{\tau}$ in each iteration. Since only components with singular values exceeding $\tau$ are retained when applying $\mathcal{D}_{\tau}$ to a matrix, an SVD package is usually called to compute  only these singular values and the corresponding singular vectors. Therefore, if the rank of matrices in each iteration is small, the algorithms can have low temporal and spatial complexity. 

There are also  some nuclear norm minimization algorithms which do not rely on the SVT, for example the Frank-Wolfe algorithm and its variants \cite{JS:ICML:10,MZWG:SISC:16}. We omit the details and interested readers should consult the references.

%% file: pgd.tex
\section{Projected gradient descent based on matrix factorization}\label{sec:pgd}
As stated in the last section, the low rank structure can be exploited effectively by nuclear norm minimization as it is amenable to detailed analysis. However, solving nuclear norm minimization by the semidefinite programming or the first order methods is computationally expensive for large scale problems. Since in an $n\times n$ rank $r$ matrix the number of degrees of freedom is $(2n-r)r$,  we can parameterize a rank $r$ matrix using a multiple of $nr$ variables. 
Alternative to convex optimization, many nonconvex algorithms  have been designed based on  the reparameterization  of low rank matrices to  solve the following variant of the rank minimization problem:
\begin{align*}
\min_{\BZ}\frac{1}{2}\|\A(\BZ)-\by\|_2^2\quad\mbox{subject to}\quad\rank(\BZ)\leq r.
\numberthis\label{eq:rank_constraint}
\end{align*}
Clearly, when $\A$ is injective on matrices of rank at most $r$, the underlying rank $r$ matrix $\BX$ is also the unique solution to \eqref{eq:rank_constraint}. In this section, we  review the nonconvex projected gradient descent (PGD) algorithm based on matrix factorization.
\begin{figure}
\centering
\includegraphics[width=10cm,height=3cm]{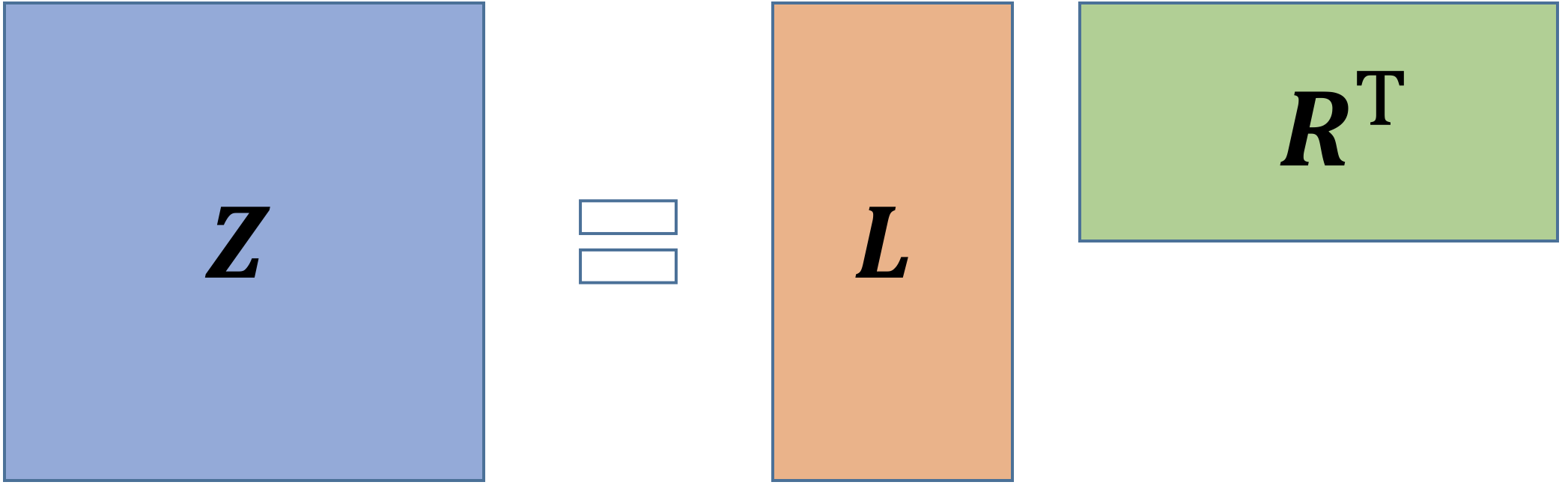}
\caption{Burer-Monteiro factorization of a low rank matrix.}\label{fig:factorization}
\end{figure}
Suppose the target rank $r$ of the underlying matrix $\BX$  is known a priori. Then it is evident that a matrix $\BZ$ has rank at most $r$ if and only if it can be factorized as a product of two rank $r$ matrices of the form (known as Burer-Monteiro factorization in optimization; see Figure~\ref{fig:factorization})
\begin{align*}
\BZ=\BL\BR^\top,\numberthis\label{eq:matrix_fact}
\end{align*}
where $\BL\in\R^{n\times r}$ and $\BR\in\R^{n\times r}$. Substituting this factorization into \eqref{eq:rank_constraint} can remove the rank constraint and turn \eqref{eq:rank_constraint} into a rank free optimization problem: 
\begin{align*}
\min_{(\BL,\BR)\in\CS} f(\BL,\BR) :=\frac{1}{2}\|\A(\BL\BR^\top)-\by\|_2^2+\bgam(\BL,\BR).\numberthis\label{eq:unconstraint1}
\end{align*}
Since the matrix factorization of the form \eqref{eq:matrix_fact} is not unique for a given matrix,
compared with the objective function in \eqref{eq:rank_constraint}, two more ingredients (i.e., a constraint set $\CS$ and a regularization function ${\bgam}(\cdot,\cdot)$) are often added in \eqref{eq:unconstraint1} to encode additional structures on the solutions we would like to seek.  Let $\BX=\BL_*\BR_*^\top$ be a {\em desired} matrix factorization of the ground truth. One typically chooses $\CS$ and $\bgam(\cdot,\cdot)$ in such a way that
\begin{align*}
(\BL_*,\BR_*)\in\CS\quad\mbox{and}\quad \bgam(\BL_*,\BR_*)=0.\numberthis\label{eq:sol_cond}
\end{align*}
Noting the fact $\A(\BL_*\BR_*^\top)=\by$, it follows that $f(\BL_*,\BR_*)=0$, so  $(\BL_*,\BR_*)$ is an optimal solution to \eqref{eq:unconstraint1}. Therefore, finding the underlying matrix $\BX$ from $\A(\BX)=\by$  can be cast as the problem of solving for the global minima of \eqref{eq:unconstraint1}. 
Moreover, a nonconvex projected gradient descent (PGD) algorithm can be developed to tackle this problem,
\begin{align*}
\begin{cases}
\widetilde{\BL}_k =  \BL_k-\alpha_k\nabla_{\BL} f(\BL_k,\BR_k)\\
\widetilde{\BR}_{k} =  \BR_k-\alpha_k\nabla_{\BR} f(\BL_k,\BR_k)\\
\lb\BL_{k+1},\BR_{k+1}\rb = \P_{\CS}\lsb(\widetilde{\BL}_k,\widetilde{\BR}_k)\rsb,
\end{cases}\numberthis\label{alg:pgd}
\end{align*}
where $\alpha_k$ is the stepsize, $\nabla_{\BL} f(\BL_k,\BR_k)$ and $\nabla_{\BR} f(\BL_k,\BR_k)$ are the partial gradients of $f$ evaluated at $(\BL_k,\BR_k)$, and $\P_{\CS}\lsb\cdot\rsb$ is the projection onto the set $\CS$.

\subsection{Recovery guarantees of PGD}
Despite the inherent nonconvex nature of \eqref{eq:unconstraint1}, theoretical recovery guarantee can be established for  PGD  with a proper initialization. A commonly used initial guess which can well approximate the underlying matrix is the so-called spectral initialization 
\begin{equation}\label{eq:specinit}
\bm{Z}_0=\mathcal{T}_r(\alpha\mathcal{A}^*(\bm{y})),
\end{equation}
where $\mathcal{A}^*$ is the adjoint of $\mathcal{A}$, $\alpha>0$ is a proper scaling factor, and $\T_r(\cdot)$ is the hard thresholding operator  which returns the best rank $r$ approximation of a given matrix (cf. the SVT in \eqref{eq:svt_operator}).  Moreover, $\T_r(\cdot)$ can be computed by the truncated SVD, 
\begin{equation}\label{eq:tsvd}
\mathcal{T}_r(\bm{Z})=\sum_{i=1}^{r}\sigma_i\bm{u}_i\bm{v}_i^T,
\quad\mbox{where }
\bm{Z}=\sum_{i=1}^{n}\sigma_i\bm{u}_i\bm{v}_i^T\mbox{ is the SVD of }\bm{Z}.
\end{equation}
We will not discuss  the approximate accuracy of  $\BZ_0$ to $\BX$ here; see \cite{WCCL:SIMAX:16,Rec:JMLR:11,CLS:TIT:15} for related results. Starting from the spectral initialization, a sufficiently close initial guess can be constructed for PGD, and then exact recovery guarantee can be established. 
\subsubsection{Matrix sensing}\label{sec:pgd_sense} When $\A$ obeys the RIP, we do not have any requirement on the unknown matrix to be reconstructed. Thus, the constraint $(\BL,\BR)\in\CS$ in \eqref{eq:unconstraint1} can be removed. Noting that $\BL\BR^\top=(\alpha\BL)(\alpha^{-1}\BR)^\top$ for any $\alpha\ne 0$, without a regularization function $\bgam(\cdot,\cdot)$, there exist solutions with $\|\BL\|_F\rightarrow \infty$ and $\|\BR\|_F\rightarrow 0$, or vice versa. This is not favorable for  the purpose of  computation and analysis. In order to avoid this situation, we can choose 
\begin{align*}
\bgam(\BL,\BR)=\lambda \|\BL^\top\BL-\BR^\top\BR\|_F^2\numberthis\label{eq:gamma}
\end{align*}
for a parameter $\lambda>0$, and then solve the following unconstraint optimization problem:
\begin{align*}
\min_{\BL,\BR}\frac{1}{2}\|\A(\BL\BR^\top)-\by\|_F^2+\lambda \|\BL^\top\BL-\BR^\top\BR\|_F^2.\numberthis\label{eq:unconstraint2}
\end{align*}
Let $\BX=\BU\BS\BV^\top$ be the compact SVD of $\BX$. Define $\BL_*=\BU\BS^{1/2}$ and $\BR_*=\BV\BS^{1/2}$. It can be easily seen that 
$\bgam(\BL_*,\BR_*)=0$, i.e., the second condition in \eqref{eq:sol_cond} is satisfied. Thus we can solve for the global minima of \eqref{eq:unconstraint2} to reconstruct $\BX$.

The gradient descent algorithm (listed in \eqref{alg:pgd} without projection) for matrix sensing is investigated in \cite{TBSSR:ICML:16} in terms of the RIP of the sensing operator. If the spectral initialization \eqref{eq:specinit} with $\alpha=1$ is used,  then the sequence  generated by \eqref{alg:pgd} converges to the global minimizer provided $\mathcal{A}$ satisfies the RIP with a small constant depending on $r$ and the condition number $\kappa$ of $\bm{X}$. To meet this condition, the number of Gaussian measurements needed is not optimal. To relax the requirement on the RIP condition so that optimal sampling complexity can be achieved, a refinement of the spectral initialization is used. The refined initialization is constructed based on $O(\log(\sqrt{r}\kappa))$ iterations of the iterative hard thresholding algorithm which will be reviewed in Section \ref{sec:IHT}. With the refined initialization, the following theorem can be established. 
\begin{theorem}[\cite{TBSSR:ICML:16}] If $\mathcal{A}$ satisfies the RIP with $\delta_{6r}< c$ for some small universal constant $c>0$, then the sequence of iterates generated by \eqref{alg:pgd} with a proper stepsize and the refined initialization converges linearly to a global minimizer which obeys $\A(\BL\BR^\top)=\by$.
\end{theorem}

\subsubsection{Matrix completion}\label{sec:pgd_completion} 
For matrix completion under uniform sampling, we are mainly interested in reconstructing a $\mu_0$-incoherent matrix, see Definition~\ref{asump:MC1}.  With  $\BL_*$ and $\BR_*$ defined in the same way as in Section~\ref{sec:pgd_sense}, the fact $\BX$ is $\mu_0$-incoherent implies that 
\begin{align*}
\|\BL_*\|_{2,\infty} \leq \sqrt{\frac{\mu_0r}{n}}\|\BX\|_2^{1/2}\quad\mbox{and}
\quad\|\BR_*\|_{2,\infty} \leq \sqrt{\frac{\mu_0r}{n}}\|\BX\|_2^{1/2}.
\end{align*}
Let $\BZ_0$ be the matrix obtained from the spectral initialization \eqref{eq:specinit} with $\alpha = n^2/m$. It can be shown that $\|\BX\|_2\leq 2\|\BZ_0\|_2$ with high probability provided $m\gtrsim \mu_0\kappa^2 r^2\log n$ \cite{ZL:ArXiv:16}. Thus, if we define 
\begin{align*}
\mathcal{C}=\left\{(\BL,\BR)~\big|~\|\bm{L}\|_{2,\infty}\leq \sqrt{\frac{2\mu_0 r}{n}}\|\bm{Z}_0\|_2^{1/2}\mbox{ and }\|\bm{R}\|_{2,\infty}\leq \sqrt{\frac{2\mu_0 r}{n}}\|\bm{Z}_0\|_2^{1/2}\right\},\numberthis\label{eq:setC}
\end{align*}
there holds $(\BL_*,\BR_*)\in\CS$, so we can choose this $\CS$ in \eqref{eq:unconstraint1}. Noting  that the unbalanced situation in Section~\ref{sec:pgd_sense} still exists here,  we can use the same $\bgam(\cdot,\cdot)$ as the regularization function.  Putting it all together, we can attempt to reconstruct the low rank factors of $\BX$ by applying PGD described in \eqref{alg:pgd} with $\CS$ and $\bgam(\cdot,\cdot)$ given in \eqref{eq:setC} and \eqref{eq:gamma} respectively. Moreover, the projection onto $\CS$ can be computed efficiently by trimming each row of $\widetilde{\BL}_k$ and $\widetilde{\BR}_k$.  That is, 
\begin{align*}
&\BL_{k+1}(i,:) = 
\begin{cases}
\widetilde{\BL}_{k}(i,:) & \mbox{if }\|\widetilde{\BL}_{k}(i,:)\|_{2,\infty}\le \sqrt{\frac{2\mu_0 r}{n}}\|\bm{Z}_0\|_2^{1/2}\\
\frac{\widetilde{\BL}_{k}(i,:) }{\|\widetilde{\BL}_{k}(i,:)\|}\sqrt{\frac{2\mu_0 r}{n}}\|\bm{Z}_0\|_2^{1/2}&\mbox{otherwise},
\end{cases}
\end{align*}
and $\BR_{k+1}$ can be computed similarly from $\widetilde{\BR}_{k}$.

Let $\BZ_0=\BU_0\BS_0\BV_0^\top$ be the compact SVD of $\BZ_0$. We can construct a provable good initial guess as follows:
\begin{align*}
\BL_0 = \P_{\CS}\lsb \BU_0\BS_0^{1/2}\rsb\quad\mbox{and}\quad\BR_0 = \P_{\CS}\lsb \BV_0\BS_0^{1/2}\rsb.\numberthis\label{eq:initial_completion}
\end{align*}
With this initial guess, the linear convergence of PGD can be established provided a sufficient number of entries are observed from the underlying matrix. 
\begin{theorem}[\cite{ZL:ArXiv:16}] Assume $\BX$ is $\mu_0$-incoherent and each pair of indices $(i_{\ell},j_{\ell})$ in $\Omega$ is sampled independently and uniformly from $\{1,\ldots,n\}\times \{1,\ldots,n\}$ with replacement.  Then with high probability the sequence of iterates generated by \eqref{alg:pgd} with a proper stepsize and the initial guess constructed by \eqref{eq:initial_completion} converges linearly to a global minimizer which obeys $\P_\Omega(\BL\BR^\top)=\P_\Omega(\BX)$ provided $m\gtrsim \mu_0\kappa^2r^2\max(\mu_0,\log n)n$.
\end{theorem}
\subsubsection{Phase retrieval}
The target matrix $\BX=\bx\bx^\top$ in phase retrieval is a rank-$1$ positive semidefinite matrix, so we can choose $\BL=\BR\in\R^{n}$ in \eqref{eq:unconstraint1}. Since $\BL$ is an $n\times 1$ vector, we replace it by the bold face lower letter $\bz$. The unbalanced situation in general matrix recovery will not appear here. In other words, $\bgam(\bz,\bz)=0$ if we choose the regularization function in \eqref{eq:gamma}. Thus without assuming any structure on $\bx$, the objective function in \eqref{eq:unconstraint1} reduces to 
\begin{align*}
f(\bz) &= \frac{1}{2}\|\A(\bz\bz^\top)-\by\|_2^2=\frac{1}{2}\sum_{\ell=1}^m(|\ba_\ell^\top\bz|^2-\by_\ell)^2.
\end{align*}
and the corresponding projected  gradient descent algorithm  can be rewritten explicitly as 
\begin{equation}\label{eq:WF}
\bm{z}_{k+1}=\bm{z}_k-\frac{\alpha_k}{m}\sum_{\ell=1}^m\left(|\bm{a}_\ell^\top\bm{z}_k|^2-y_\ell\right)(\bm{a}_\ell^\top\bm{z}_k)\bm{a}_\ell,
\end{equation}
where $\alpha_k$ is the stepsize. In the complex case, the gradient should be calculated using Wirtinger calculus, so the gradient descent iteration is also referred to as Wirtinger flow in the literature \cite{CLS:TIT:15}. For Wirtinger flow, the initial guess can also be constructed from the spectral initialization in \eqref{eq:specinit} with $r=1$ and $\alpha=1$:  Let $\BZ_0=\bz_0\bz_0^\top$ and then rescale $\bz_0$ such that $\|\bz_0\|_2=\|\by\|_1/m$. With this initialization, the theoretical guarantee of Wirtinger flow can be established. 
\begin{theorem}[\cite{CLS:TIT:15}]
Assume $\ba_{\ell}\sim\N(0,\BI_n)$ and $\by=|\BA\bx|^2$. Then with high probability Wirtinger flow with a proper stepsize and the initial guess constructed from the spectral initialization converges linearly to $\bx$ provided $m\gtrsim n\log n$.
\end{theorem}

\paragraph{Remark}
For conciseness, we have discussed the simplest  PGD algorithm in this section, and yet many other algorithms can be developed based on the matrix factorization model, for example alternating minimization \cite{JNS:TOC:13,Har:FOCS:14,wyz2012lmafit} and alternating steepest descent  \cite{tannerwei2016asd}.  In particular, a large family of related algorithms have been discussed in \cite{SL:TIT:16}. For phase retrieval, there have been many variants of Wirtinger flow with improved computational efficiency or sampling complexity \cite{CC:CPAM:17,WGE:TIT:18,ZCL:ICML:16}. For example, a truncated variant of Wirtinger flow based on Poisson loss was shown to be able to converge to $\bx$ linearly provided  $m\gtrsim n$ and a truncated initialization is used \cite{CC:CPAM:17}.  In addition, if we utilize matrix factorization with three blocks,  Grassmann manifold algorithms can be developed for low rank matrix recovery \cite{KMO:TIT:10,NS:NIPS:12,BA:NIPS:11,mishra2012riemannian,mishra2014fixed,mishra2014r3mc}.

Exact recovery guarantees have been presented for PGD with a proper initialization. Inspired by the observation that PGD seeded with a random guess often converges to a global minimizer, another line of research has been devoted to study the geometric landscape of the objective function $f$ in \eqref{eq:unconstraint1} \cite{GJZ:ICML:17,SQW:FCM:18,GJZ:ICML:17,GJZ:ArXiv:17}. Typical results are $f$ does not have a spurious local minima and there exist a descent direction at each saddle point, so that any algorithm which can converge to a local minimizer is able to find a global minimizer. Moreover, many algorithms have 
been designed to escape saddle points efficiently \cite{GHJY:COLT:15,JGNKJ:ArXiv:17,CDHS16,AAZBHM16}.

%% file: manifold.tex
\section{Algorithms on embedded manifold of low rank matrices}\label{sec:manifold}
We have already seen that matrix factorization and the corresponding nonconvex algorithms can be utilized to exploit the structure in low rank matrix recovery effectively and efficiently. In this section,  another class of nonconvex algorithms to exploit the low rank structure are presented, which  
proceed by minimizing a smooth loss function over the embedded manifold low rank matrices, 
\begin{equation}\label{eq:LScons}
\min_{\bm{Z}\in\mathcal{M}_r}~\|\mathcal{A}(\bm{Z})-\bm{y}\|_2^2.
\end{equation}
Here $\M_r$ denotes the set of fixed rank $r$ matrices. It is well-known that $\M_r$ is a smooth manifold \cite{Van:SIOPT:13}. We begin our discussion with the simple iterative hard thresholding algorithm for \eqref{eq:LScons} and then extend it to a class of Riemannian optimization algorithms.

\subsection{Iterative hard thresholding}\label{sec:IHT}
The objective function in \eqref{eq:LScons} is convex and smooth. Although the set $\mathcal{M}_r$ is non-convex, the projection onto it has a closed form and can be computed by the truncated SVD; see \eqref{eq:tsvd}. 
Thus, a simple algorithm for \eqref{eq:LScons} is the following iterative hard thresholding (IHT) algorithm:
\begin{equation}\label{eq:IHT}
\begin{cases}
\BG_k=\mathcal{A}^\top\left(\bm{y-\mathcal{A}(\bm{Z}_k)}\right)\\
\bm{Z}_{k+1}=\mathcal{T}_{r}\left(\bm{Z}_k-\alpha_k\BG_k\right),
\end{cases}
\end{equation}
where $\alpha_k$ is the stepsize. In each iteration, IHT first computes the gradient descent direction $\BG_k$ of the quadratic  objective function and then updates the current estimate $\BZ_k$ along $\BG_k$, followed by projection onto $\M_r$ via the hard thresholding operator $\T_r$. IHT was first designed for compressed sensing in \cite{bludav2009iht} and then extended to low rank matrix recovery in  \cite{JMD:NIPS:10} (referred to as SVP in there).
Theoretical recovery guarantee of IHT  was established in  \cite{JMD:NIPS:10} in terms of the RIP of $\A$, showing that IHT is able to reconstruct a rank-$r$ matrix provided that $\mathcal{A}$ satisfies the RIP with the constant $\delta_{2r}<1/3$ and  the stepsize is chosen to be $\alpha_k=1/(1+\delta_{2k})$.

We can also choose the search stepsize in an adaptive way. Since the objective function is a least-squares, an exact line search in a linear subspace leads to a stepsize with a closed form. In particular, it is proposed  in \cite{TW:SISC:13} to do exact line search in the column subspace of $\BZ_k$: $\{\bm{U}_k\bm{B}^\top|\bm{B}\in\mathbb{R}^{n\times r}\}$, where $\bm{U}_k\in\mathbb{R}^{n\times r}$ consists of  the left $r$ singular vectors of $\bm{Z}_k$. Noting that $\BG_k$ is the gradient descent direction, the stepsize for the exact line search along  the projection of $\BG_k$ onto the column subspace is given by
\begin{equation}\label{eq:IHT_stepsize}
\alpha_k=\frac{\|\bm{U_k}\bm{U}_k^T\bm{G}_k\|_F^2}{\|\mathcal{A}(\bm{U_k}\bm{U}_k^T\bm{G}_k)\|_2^2}.
\end{equation}
Other subspaces such as the row subspace of $\BZ_k$ can also be used to compute the stepsize. The algorithm \eqref{eq:IHT} with the adaptive stepsize is known as normalized iterative hard thresholding (NIHT). It is proven in \cite{TW:SISC:13} that, if $\mathcal{A}$ satisfies the  RIP with the constant $\delta_{3r}<1/5$, NIHT converges linearly to $\bm{X}$, which is optimal in sampling complexity under Gaussian measurements. The result  in \cite{TW:SISC:13} applies equally for a constant stepsize and thus does not rely on some unknown stepsize in contrast to the one in \cite{JMD:NIPS:10}.

Despite the optimal recovery guarantee of SVP and NIHT, they suffer from the slow asymptotic convergence rate of gradient descent methods. To improve the efficiency, one may consider conjugate gradient descent type methods. A family of conjugate gradient iterative hard thresholding (CGIHT) algorithms were proposed in \cite{BTW:II:15}.  
It was also proved  that a restarted version of CGIHT converges linearly to $\bm{X}$ under the RIP assumption of $\mathcal{A}$.

The performance guarantee of IHT for matrix completion is recently investigated in \cite{DingChen18} using the leave-one-out analysis. To the best of our knowledge, IHT for phase retrieval has not been studied yet. We will omit further details of IHT because in each iteration the SVD on an $n\times n$ matrix is needed to compute the projection onto $\M_r$ which is computationally inefficient. Next, we will see how to modify IHT in an elegant way to improve the computational efficient which leads to a class of Riemannian optimization algorithms.

\subsection{Riemannian optimization on low rank manifold}
We first  refer the reader to the textbook by \cite{AbMaSe2008manifold} for  comprehensive treatments of Riemannian optimization. Here we investigate a Riemannian optimization algorithm for low rank matrix recovery based on $\M_r$, which is a smooth Riemannian manifold when  embedded into the Euclidean space $\mathbb{R}^{n\times n}$ with the standard inner product. A Riemannian conjugate gradient descent algorithm was first introduced into matrix completion in \cite{Van:SIOPT:13}. The difference and connection between the Riemannian optimization on the embedded manifold of fixed rank $r$ matrices and IHT were pointed out in \cite{KeWeiThesis}, and then exact recovery guarantees of  the corresponding Riemannian optimization algorithms were established in \cite{WCCL:SIMAX:16,WCCL:ArXiv:16,CaiWeiphase} for matrix sensing, matrix completion and phase retrieval respectively based on the connection with IHT. 

\begin{figure}[ht!]\centering
\includegraphics[height=5cm,width=6cm]{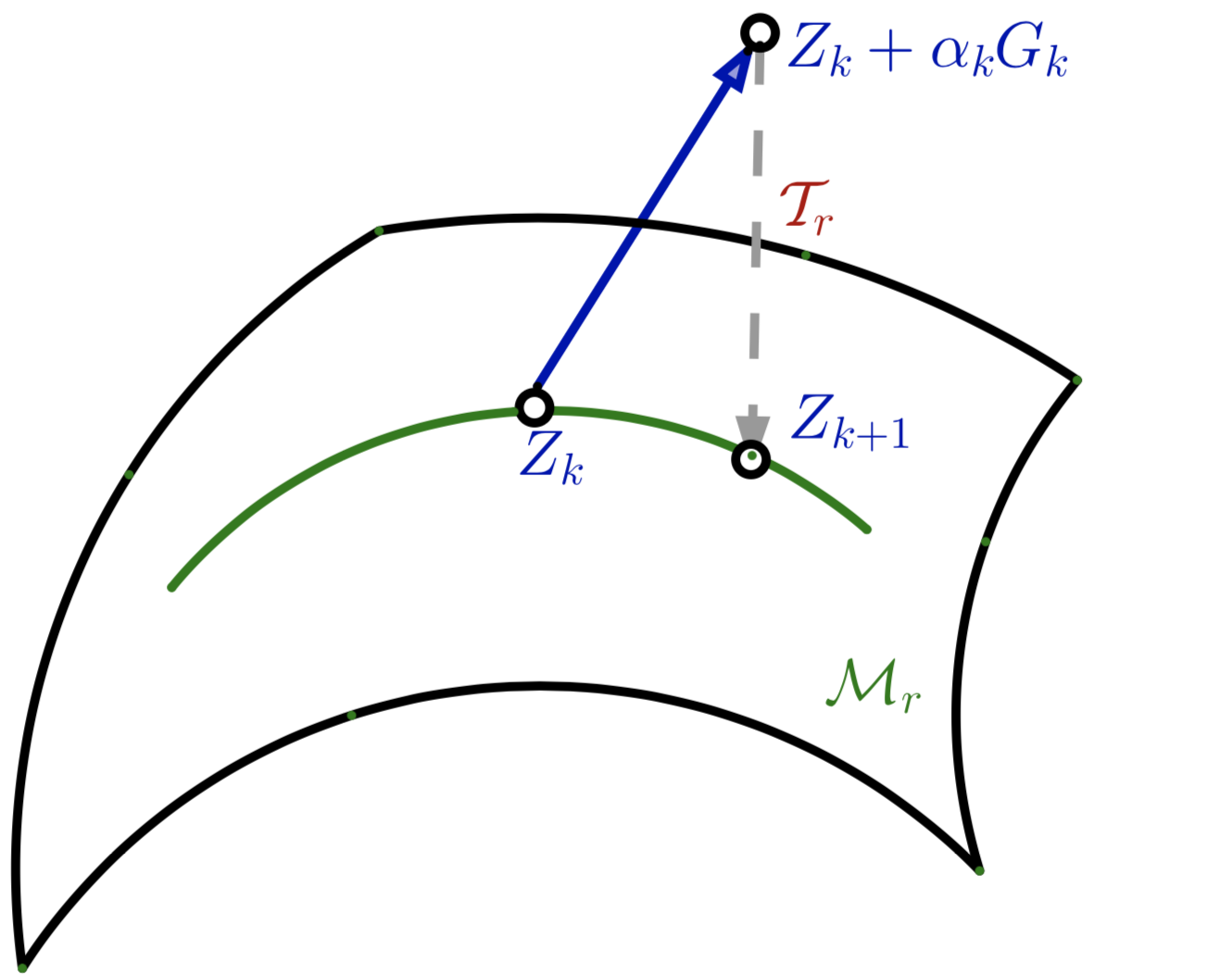}\hspace{1cm}
\includegraphics[height=5cm,width=6.5cm]{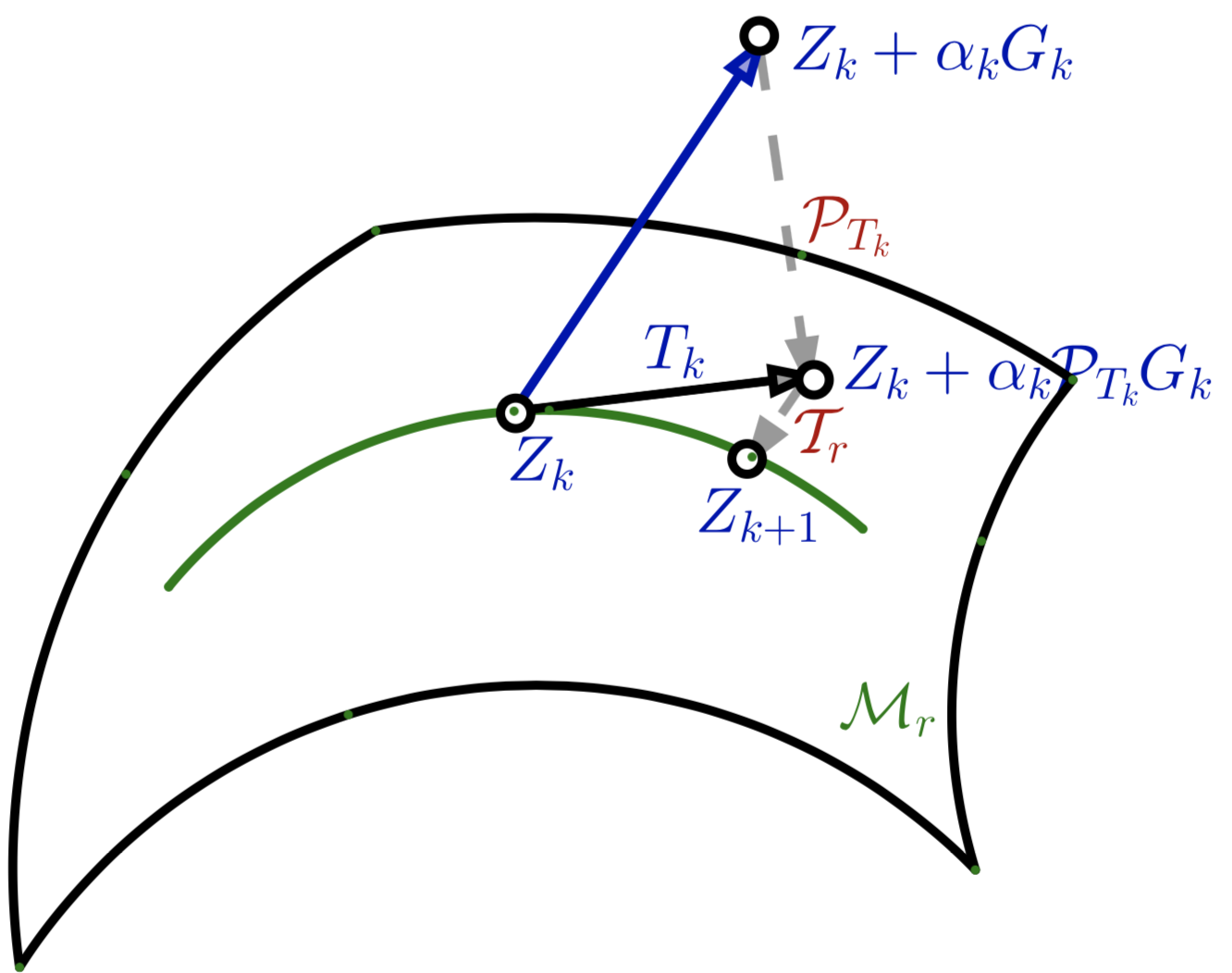}
\caption{Pictorial illustration of IHT (left) and RGrad (right).}\label{fig:ihtVSrgrad}
\end{figure}
In each iteration of IHT, we need to compute the SVD of an $n\times n$ matrix and the computational cost is $O(n^3)$ in general as  the  matrix  after the gradient descent update is typically unstructured. To overcome the high computational cost of the SVD, we can first project the matrix obtained after the gradient descent onto a low dimensional subspace, followed by projection onto the low rank matrix manifold $\M_r$. After the projection onto a low dimensional subspace, it is possible that the resulting matrix will be low rank and structured so that the projection onto $\M_r$ by the SVD can be computed efficiently. If the low dimensional subspace is selected to be the tangent space of the manifold $\M_r$ at the current estimate, we obtain the Riemannian gradient descent algorithm which is referred to as {RGrad} in the survey. The algorithm can be formally described as follows: 
\begin{equation}\label{eq:RGD}
\begin{cases}
\BG_k=\mathcal{A}^\top\left(\bm{y-\mathcal{A}(\bm{Z}_k)}\right)\\
\bm{Z}_{k+1}=\mathcal{T}_{r}\left(\bm{Z}_k+\alpha_k\P_{T_k}(\BG_k\right)),
\end{cases}
\end{equation}
where $\alpha_k$ is the stepsize,  $T_k$ is tangent space of $\M_r$ at $\BZ_k$, $\P_{T_k}$ is the associated projection operator, and in the second line, 
\begin{align*}
\P_{T_k}\lb\bm{Z}_k+\alpha_k\BG_k\rb = \bm{Z}_k+\alpha_k\P_{T_k}(\BG_k)
\end{align*}
as we will see $\BZ_k\in T_k$. In Riemannian optimization, $\T_r$ is known as a type of {\em retraction}; see \cite{AbMaSe2008manifold} for other choices of retractions. In RGrad, an exact linear search along $\P_{T_k}(\BG_k)$ yields a closed form stepsize given by 
\begin{align*}
\alpha_k = \frac{\|\P_{T_k}(\BG_k)\|_F^2}{\|\A(\P_{T_k}(\BG_k))\|_2^2}.\numberthis\label{eq:rgrad_step}
\end{align*}

 Compared with IHT, there is an additional projection onto the tangent space $T_k$ in RGrad; see Figure~\ref{fig:ihtVSrgrad} for all illustration. Due to this subtle difference, the computational efficiency can be improved significantly. 
Let $\BZ_k=\BU_k\BS_k\BV_k^\top$ be the SVD of $\BZ_k$. The tangent space $T_k$ is given by \cite{Van:SIOPT:13}
$$
T_k=\left\{\bm{U}_k\bm{B}^\top+\bm{C}\bm{V}_k^\top~|~\bm{B},\bm{C}\in\mathbb{R}^{n\times r}\right\}.
$$
It follows immediately that $\BZ_k\in T_k$. Moreover, each matrix $\BW_k=\bm{U}_k\bm{B}^\top+\bm{C}\bm{V}_k^\top$ in $T_k$ is rank at most $2r$ and a simple algebra yields 
\begin{align*}
\bm{W}_k=\begin{bmatrix}
\BU_k & \BQ_2
\end{bmatrix}
\begin{bmatrix}
\BM & \BR_1^\top\\
\BR_2 & \bm{0}
\end{bmatrix}
\begin{bmatrix}
\BV_k & \BQ_1
\end{bmatrix}^\top
\end{align*}
for matrices $\BM\in\R^{r\times r},~\BR_1\in\R^{r\times r},~\BR_2\in\R^{r\times r}$, $\BQ_2\in\R^{n\times r}$ obeying $\BQ_2\perp\BU_k$, and $\BQ_1\in\R^{n\times r}$ obeying $\BQ_1\perp\BV_k$, all of which  can be computed from $\BU_k$, $\BV_k$, $\BB$ and $\BC$ using a few matrix products. Thus, both $\begin{bmatrix}
\BU_k & \BQ_2
\end{bmatrix}$ and $ \begin{bmatrix}
\BV_k & \BQ_1
\end{bmatrix}$  are orthogonal matrices and the SVD of $\BW_k$ can be computed efficiently from 
the SVD of the middle $2r\times 2r$ matrix. The total computational cost of the SVD is $O(nr^2+r^3)$ flops, which is much smaller than $O(n^3)$ when $r\ll n$; see \cite{WCCL:SIMAX:16,Van:SIOPT:13} for details.

In addition, one can easily modify RGrad to have the Riemannian conjugate gradient descent  algorithm: 
\begin{align*}
\BZ_{k+1} = \T_r\lb\BZ_k+\P_{T_k}(\BP_k)\rb,
\end{align*}
where the new search $\BP_k$ is a weighted sum of the gradient descent direction $\BG_k$ and the previous search direction $\BP_{k-1}$. Several choices of the combination weight are available \cite{WCCL:SIMAX:16,Van:SIOPT:13}. In each iteration,  the Riemannian conjugate gradient descent algorithm has the same dominant computational cost as RGrad but with substantially faster convergence rate. The details will be omitted here. 
\subsection{Recovery guarantees of RGrad}
In this section, we present the recovery guarantees of RGrad for matrix sensing, matrix completion and phase retrieval. 

\paragraph{Matrix sensing} It was shown in \cite{WCCL:SIMAX:16} that, if $\mathcal{A}$ satisfies the RIP with $\delta_{3r}<1/(12\kappa\sqrt{r})$, then RGrad with  the spectral initialization converges linearly to $\bm{X}$. For Gaussian measurements, it implies 
$m\gtrsim \kappa^2nr^2$ sampling complexity which is suboptimal. To remedy this problem, we can follow the approach in \cite{TBSSR:ICML:16} and run  $O(\log r)$ iterations  of IHT to construct a more accurate initial guess.  Then the sampling complexity will be optimal.
\begin{theorem}[\cite{WCCL:SIMAX:16}] If $\mathcal{A}$ satisfies the RIP with $\delta_{3r}< c$ for some small absolute numerical constant $c>0$, then the sequence of iterates generated by \eqref{eq:RGD} with a proper stepsize  and the refined initialization converges linearly to $\BX$.
\end{theorem}

\paragraph{Matrix completion} Under the assumptions that $\BX$ is $\mu_0$-incoherent and the indices for the observed entries are sampled independently and uniformly with replacement, it was shown in \cite{WCCL:ArXiv:16} that RGrad with the spectral initialization converges linearly to $\BX$ with high probability provided $m\gtrsim \mu_0\kappa n^{1.5}r\log^{1.5} n$. The sampling complexity is undesirable with $n$. In order to improve this result, a refined initialization is proposed in \cite{WCCL:ArXiv:16}  which runs RGrad one pass on $O(\log n)$ nonoverlapping partitions of the observed entries followed by trimming. The following theorem can be established with the refined initialization.
\begin{theorem}[\cite{WCCL:ArXiv:16}] Assume $\BX$ is $\mu_0$-incoherent and each pair of indices $(i_{\ell},j_{\ell})$ in $\Omega$ is sampled independently and uniformly from $\{1,\ldots,n\}\times \{1,\ldots,n\}$ with replacement.  Then with high probability the sequence of iterates generated by \eqref{eq:RGD} with a proper stepsize and the refined initialization converges linearly to $\BX$ provided $m\gtrsim \mu_0\kappa^6r^2n\log^2n$.
\end{theorem}

\paragraph{Phase retrieval} Recall that the target matrix $\BX=\bx\bx^\top$ in phase retrieval is a rank-$1$ positive semidefinite matrix. RGrad in \eqref{eq:RGD} can preserve this structure in each iteration. Assume $\BZ_k$ is a rank-$1$ positive semidefinite matrix in the $k$-th iteration. Then it has the following eigenvalue decomposition
\begin{align*}
\BZ_k=\sigma_k\bu_k\bu_k^\top,
\end{align*}
where $\sigma_k\geq 0$ and $\bu_k$ is a unit vector. The tangent space of positive rank-$1$ matrices at $\BZ_k$ is given by \cite{HuangGaZh16}
\begin{align*}
T_k =\{\bu_k\bb^\top+\bb\bu_k^\top~|~\bb\in\R^n\}.
\end{align*}
Noting the special property of $T_k$, after updating $\BZ_k$ along the direction $\P_{T_k}(\BG_k)$, we can compute the new estimate $\BZ_{k+1}$ as the best rank-$1$ positive semidefinite approximation via the eigenvalue decomposition. 

Under the  Gaussian sampling model, the measurement matrices $\bm{A}_\ell=\bm{a}_\ell\bm{a}_\ell^\top$ is an outer product of two Gaussian vectors, so $\|\mathcal{A}(\bm{Z})\|_2^2$ contains the $4$-th moment of Gaussian random variables which does not possess a good concentration around its expectation. Therefore, it is not very clear how to establish the convergence of RGrad. Despite this, a truncated variant of RGrad with competitive performance was proposed in \cite{CaiWeiphase} which was able to achieve exact recovery with high probability based on the Gaussian measurement model.
\begin{theorem}[\cite{CaiWeiphase}]
Assume $\ba_{\ell}\sim\N(0,\BI_n)$ and $\by=|\BA\bx|^2$. Then with high probability a truncated variant of RGrad with a proper stepsize and initial guess converges linearly to $\bx$ provided $m\gtrsim n$.
\end{theorem}

\subsection{PGD vs RGrad: An illustration on matrix completion}
Overall, PGD and RGrad have similar per iteration computational cost, so they are two  equally effective ways to exploit the low rank structure in low rank matrix recovery. We consider matrix completion as an illustration. The dominant per iteration computational cost of PGD  for \eqref{eq:unconstraint1} with $\CS$ in \eqref{eq:setC} and $\bgam(\cdot,\cdot)$ in \eqref{eq:gamma} is 
$O(|\Omega|r+|\Omega|+nr^2+nr)$ \cite{ZL:ArXiv:16} while that of GRrad for \eqref{eq:LScons} is $O(|\Omega|r+|\Omega|+nr^2+nr+r^3)$ \cite{WCCL:ArXiv:16}, where $|\Omega|$ denotes the number of observed entries.   

We  evaluate the performance of PGD and RGrad via a set of simple experiments. As suggested by \cite{ZL:ArXiv:16},  the regularization function and the projection are not included when implementing PGD  since the algorithm works equally well without those two components. The stepsize in PGD is determined via backtracking while the stepsize in RGrad is computed via \eqref{eq:rgrad_step}. The initial guesses are constructed from the spectral initialization \eqref{eq:specinit} with $\alpha=n^2/m$ for both algorithms.  The experiments are conducted on a
Mac Pro laptop with 2.5GHz quad-core Intel Core i7 CPUs and 16 GB memory and executed from Matlab 2014b.

\begin{figure*}[ht!]
\centering
\includegraphics[width=3in, trim = 0 7cm 0 7cm]{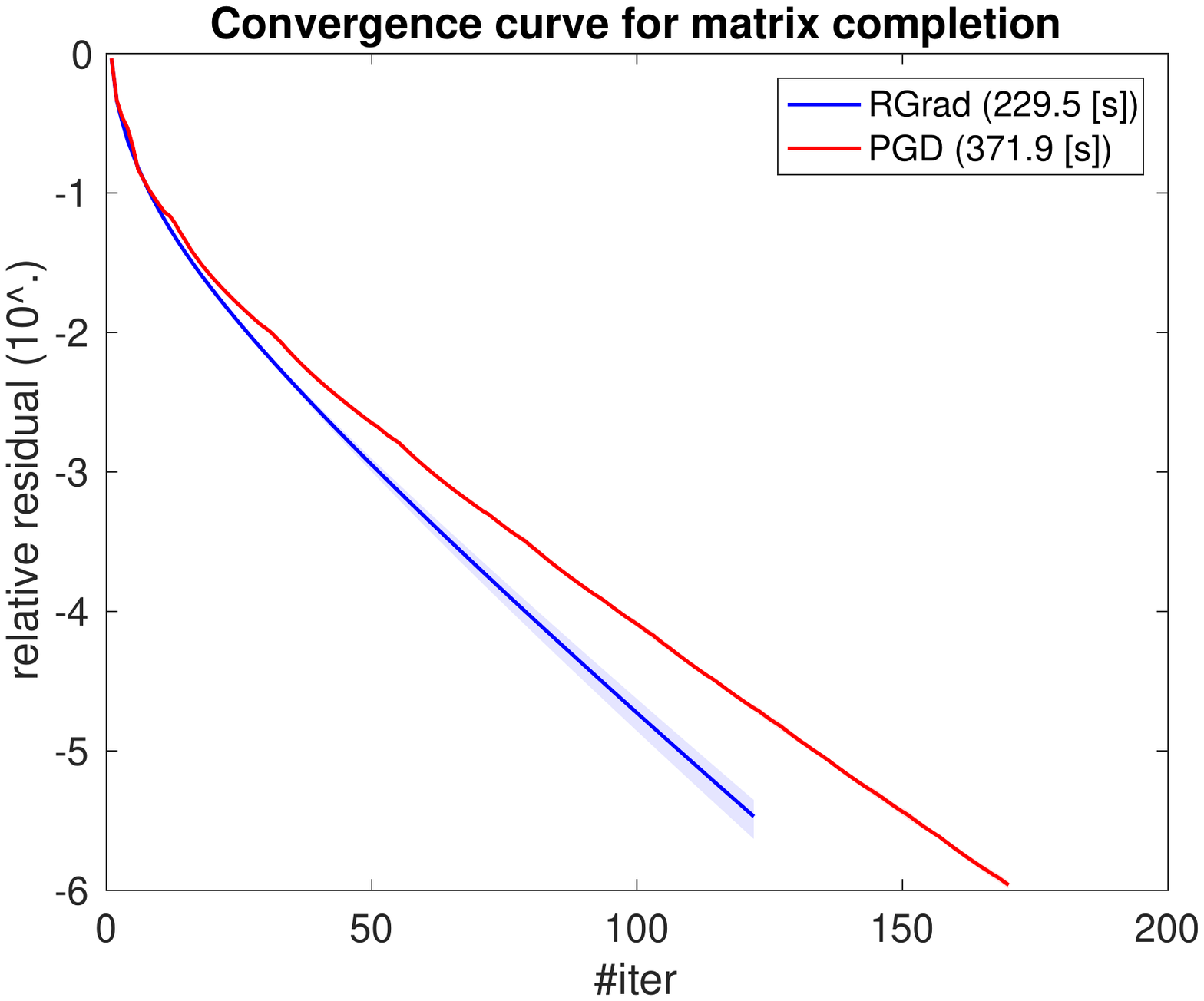}
\includegraphics[width=3in, trim = 0 7cm 0 7cm]{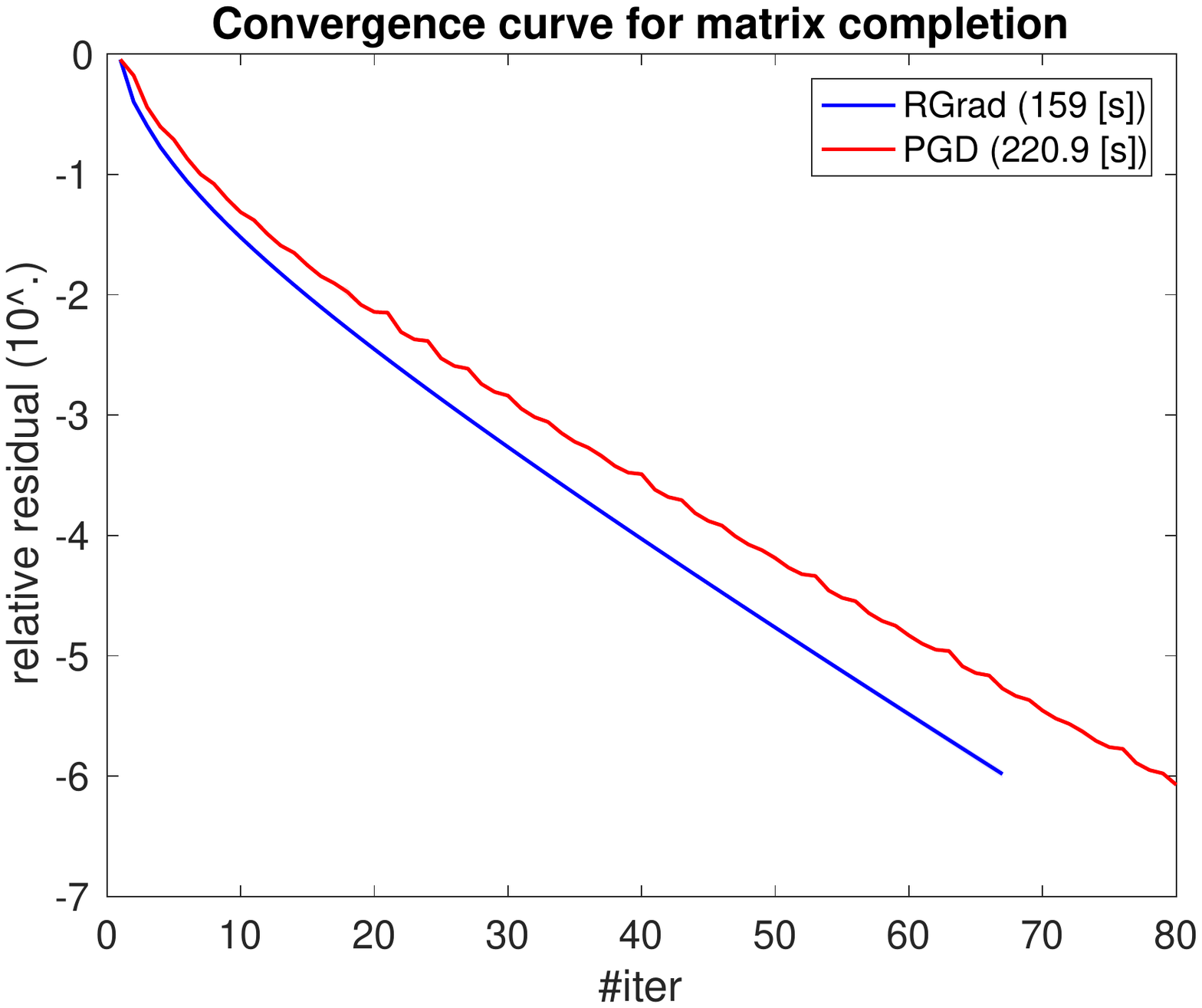}
\caption{{Relative residual (mean and standard deviation over ten random tests) as function of
number of iterations for $n=8000$, $r=100$, $m=2(2n-r)r$ (left) and $m=3(2n-r)r$ (right). The values after each algorithm are the average computational time
(seconds) for convergence.}}
\label{fig_time}
\end{figure*}

We test the algorithms on randomly generated matrices of size $8000\times 8000$ 
and rank $100$, which are computed via $\BX=\BL\BR^\top$ with $\BL$ and $\BR$ having i.i.d Gaussian entries. Two values of $m$: $m=2(2n-r)r$ and $m=3(2n-r)r$ are tested and 
the algorithms are terminated when the relative residual is less than $10^{-6}$. The relative residual plotted against
the number of iterations and the average recovery time are presented in Figure~\ref{fig_time}. It can be observed that in the setting of our tests RGrad is slightly faster than PGD, but overall they exhibit similar convergence behavior. It is worth noting that the Riemannian conjugate gradient descent algorithm whose convergence curve is not presented in the figure can be significantly faster than PGD and RGrad.

\subsection{Extensions}
Note that the key difference between IHT and RGrad is the additional projection onto a low dimensional subspace before the projection onto the low rank matrix manifold. This idea turns out to be very useful in designing fast algorithms for more general low rank matrix recovery problems. In this subsection, we give two more examples.
\subsubsection{Spectrally sparse signal reconstruction}
In many applications, the signal of interest is not low rank itself, but will exhibit a low rank structure after some linear or nonlinear transforms. A typical example is the spectrally sparse signal which appears in a wide range of applications, including magnetic resonance imaging \cite{MRI}, fluorescence microscopy \cite{Microscopy}, radar imaging \cite{Radar}, nuclear magnetic resonance (NMR) spectroscopy \cite{QMCCO:ACIE:15}. In the simplest one dimensional case, a spectrally sparse signal  
$\bx\in\C^n$ is in the form of 
\begin{align*}
\bx = \begin{bmatrix}
x_0\\
x_1\\
x_2\\
\vdots\\
x_{n-2}\\
x_{n-1}
\end{bmatrix}= \begin{bmatrix}
w_1^0 & w_2^0 & \cdots & w_r^0\\
w_1^1 & w_2^1 & \cdots & w_r^1\\
w_1^2 & w_2^2 & \cdots & w_r^2\\
\vdots & \vdots & \vdots & \vdots\\
w_1^{n-2} & w_2^{n-2} & \cdots & w_r^{n-2}\\
w_1^{n-1} & w_2^{n-1} & \cdots & w_r^{n-1}
\end{bmatrix}
\begin{bmatrix}
d_1\\
d_2\\
\vdots\\
d_r
\end{bmatrix},\numberthis\label{eq:spectrally_sparse_signal}
\end{align*}
where $w_j = e^{2\pi if_j -\tau_j},~j=1,\cdots, r$ for $r$ distinct frequencies $f_j\in [0, 1)$  and  $r$ real damping factors $\tau_j\geq 0$.

Spectrally sparse signal reconstruction or spectral compressed sensing is about reconstructing a spectrally sparse signal from the partial observed entries of the signal. Let $\Omega$ be a subset of $\{0,\cdots,n-1\}$ corresponding to the observed entries, and let $\P_{\Omega}$ be the associated sampling operator. Then the goal is to reconstruct $\bx$ from $\P_\Omega(\bx)$. In general, this is an ill-posed problem as one can fill in any values into the locations of the unknown entries. However, there is a low rank structure hidden in $\bx$ which can be utilized to complete the reconstruction task. 

\begin{figure}
\centering
\includegraphics[height = 5cm,width=6cm]{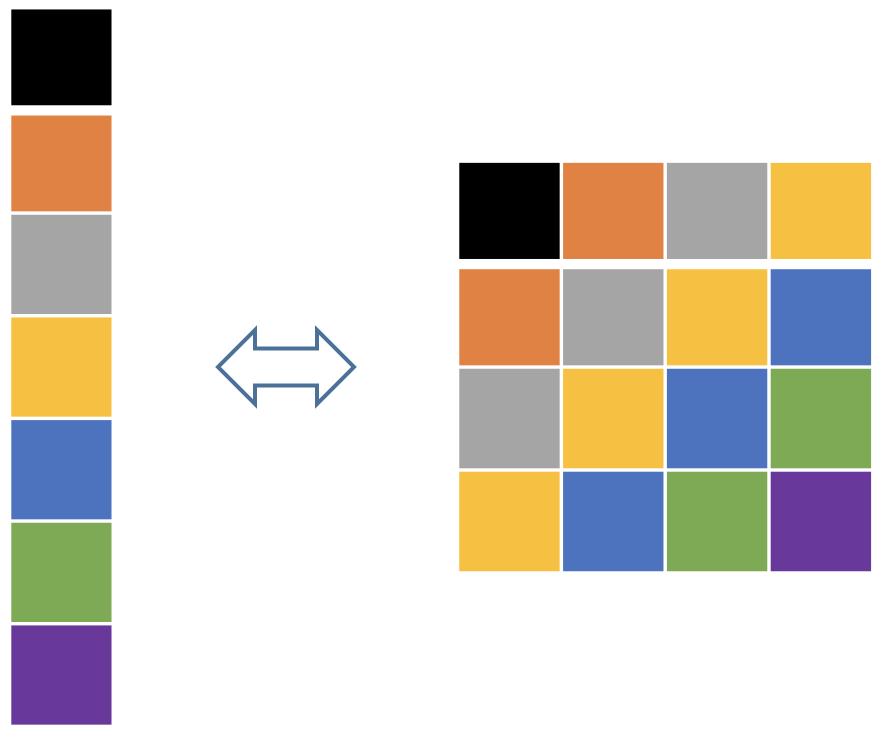} 
\caption{Vector $\bz$ (left) and its Hankel transform $\H\bz$: $\H$ maps each entry of $\bz$ into an anti-diagonal of $\H\bz$. Thus, there is a one to one correspondence between the entries of $\bz$ and the anti-diagonals of $\H\bz$.}\label{fig:hankel}
\end{figure}
\begin{figure}[ht!]\centering
\includegraphics[height=5cm,width=6.5cm]{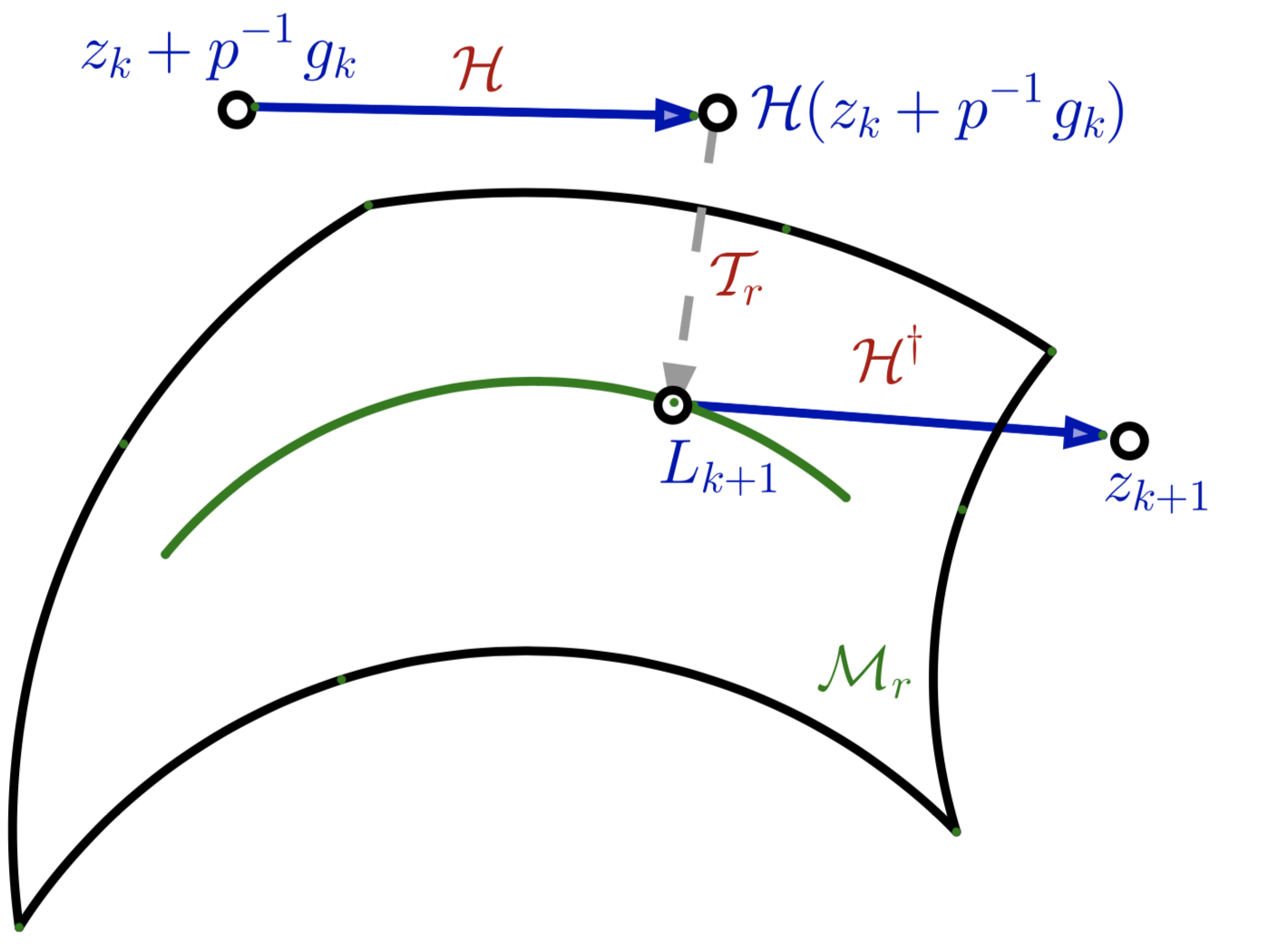}\hspace{1cm}
\includegraphics[height=5cm,width=6.5cm]{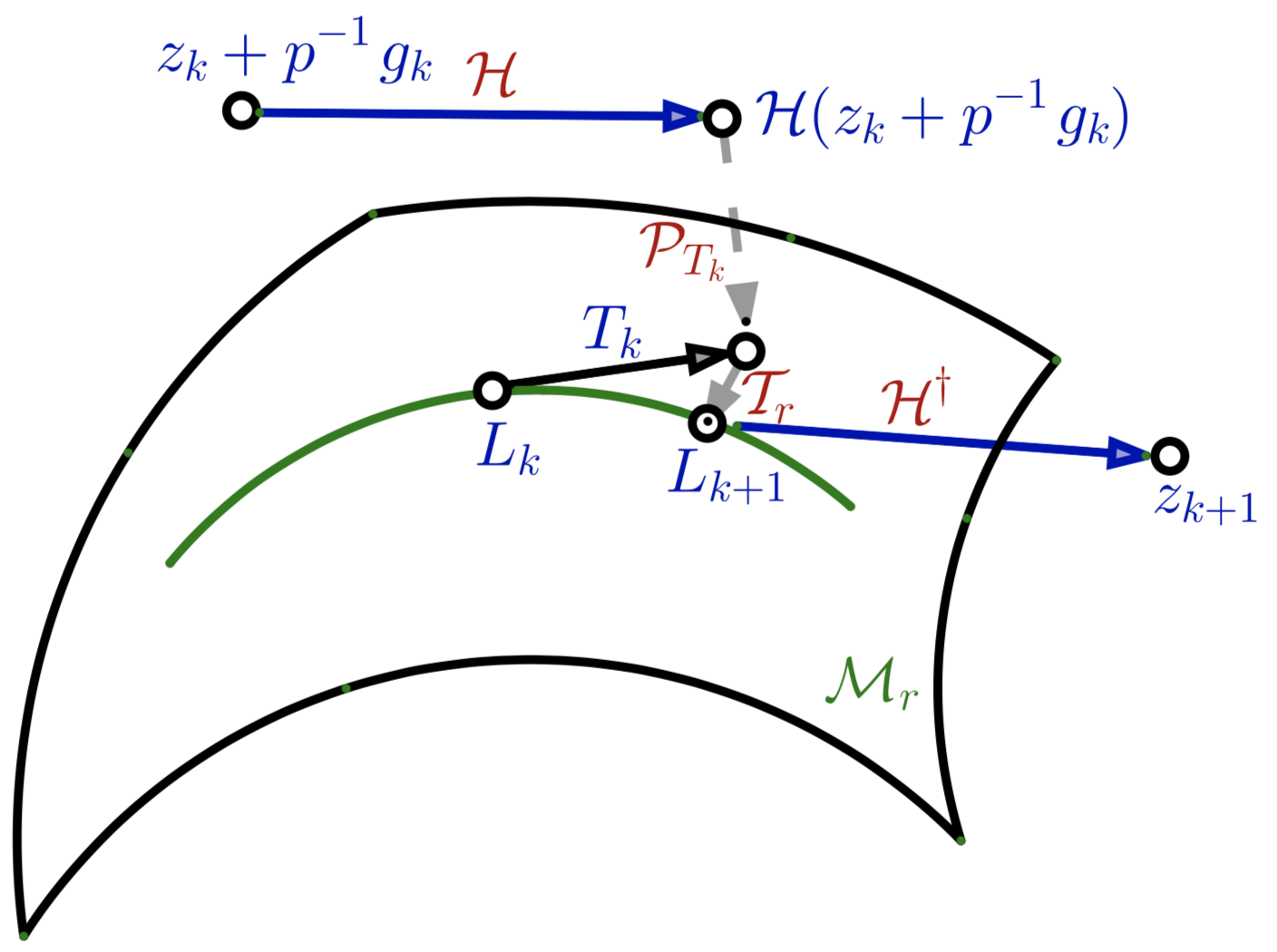}
\caption{Pictorial illustration of IHT (left) and FIHT (right) for spectrally sparse signal reconstruction.}\label{fig:ihtVSfiht}
\end{figure}

Given a vector $\bz\in\C^n$, let $\H$ be a linear operator which maps $\bz$ into an $n_1\times n_2$ Hankel  matrix obeying $n_1+n_2=n+1$ (see Figure~\ref{fig:hankel}), 
$$
[\mathcal{H}\bm{z}]_{ij}=z_{i+j}, \quad \forall ~i\in\{0,\ldots,n_1-1\},~j\in\{0,\ldots,n_2-1\}.
$$
Because $\bx$ is a spectrally sparse signal, a simple calculation can show that 
$\H\bx$ admits the following Vandermonde decomposition:
\begin{align*}
\H\bx  = \begin{bmatrix}
1 & 1 & \cdots & 1\\
w_1 & w_2 & \cdots & w_r\\
\vdots&\vdots &\vdots &\vdots\\
w_1^{n_1-1} & w_2^{n_1-1} &\cdots &w_r^{n_1-1}
\end{bmatrix}
\begin{bmatrix}
d_1 & & &\\
& d_2 & & \\
& & \ddots & \\
& & & d_r
\end{bmatrix}
 \begin{bmatrix}
1 & 1 & \cdots & 1\\
w_1 & w_2 & \cdots & w_r\\
\vdots&\vdots &\vdots &\vdots\\
w_1^{n_2-1} & w_2^{n_2-1} &\cdots &w_r^{n_2-1}
\end{bmatrix}.
\end{align*}
From this decomposition, one can easily see that $\rank(\H\bx)= r$, so $\H\bx$ is a low rank matrix when $r\ll n_1$ and $r\ll n_2$. Thus we can attempt to reconstruct $\bx$ by seeking a signal which  fits the observed entries as well as possible and at the same time is low rank after Hankel transform:
\begin{align*}
\min_{\bm{z}}\|\P_{\Omega}(\bm{z})-\P_{\Omega}(\bm{x})\|_2^2
\quad\mbox{subject to}\quad
\mathrm{rank}(\mathcal{H}\bm{z})=r.\numberthis\label{eq:hankel_prob}
\end{align*}
There is no closed-form projection onto the feasible set $\{\bm{z}~|~\mathrm{rank}(\mathcal{H}\bm{z})=r\}$, so  projected gradient descent is not directly applicable.
In \cite{CWW:ACHA:18}, an approximate projected gradient descent algorithm, still referred to as IHT, is proposed for \eqref{eq:hankel_prob}:
\begin{equation}\label{eq:IHTHankel}
\begin{cases}
 \bm{g}_k=\P_\Omega(\bm{x}-\bm{z}_{k})\\
 \bm{z}_{k+1}=\H^\dag\T_r\H(\bm{z}_{k}+\alpha_k\bm{g}_k),
\end{cases}
\end{equation}
where $\alpha_k$ is the stepsize and $\mathcal{H}^{\dag}$ is the pseudo-inverse of $\mathcal{H}$. In each iteration, IHT first updates the current estimate $\bz_k$ along the gradient descent direction $\bg_k$. Then  the Hankel matrix corresponding to
the update is formed via the application of the Hankel transform $\H$, followed by the SVD truncation
to the best rank-$r$ approximation via the hard thresholding operator $\T_r$. Finally, 
the new estimate $\bz_{k+1}$ is obtained via the application of pseudo-inverse Hankel transform $\H^\dagger$. See Figure~\ref{fig:ihtVSfiht} (left) for an illustration.

In order to reduce the computational cost of the SVD in IHT, inspired by RGrad, we can  add an additional subspace projection before truncating the Hankel matrix to its nearest rank-$r$ approximation. This leads to the FIHT algorithm proposed in \cite{CWW:ACHA:18}: 
\begin{equation}\label{eq:IHTHankel}
\begin{cases}
 \bm{g}_k=\P_\Omega(\bm{x}-\bm{z}_{k})\\
 \bm{z}_{k+1}=\H^\dag\P_{T_k}\T_r\H(\bm{z}_{k}+\alpha_k\bm{g}_k),
\end{cases}
\end{equation}
where $T_k$ is selected to be the tangent space of the rank $r$ matrix manifold $\M_r$ at the previous rank-$r$ matrix $\BL_k$; see Figure~\ref{fig:ihtVSfiht} (right).

As in RGrad, the truncation to the rank-$r$ matrix manifold $\M_r$ in FIHT can be computed very efficiently. Thus, FIHT is computationally much faster than IHT. For example, numerical simulation shows FIHT can reconstruct a $128\times 128\times 1024$ three dimensional spectrally sparse signal with $20$ frequencies from $4\%$ of the known entries in less than an hour on a laptop \cite{CWW:ACHA:18}. Moreover, exact recovery guarantee of FIHT can also be established, which shows under the sampling with replacement model FIHT with a proper initialization  can achieve successful recovery with high probability provided $\H\bx$ is well conditioned and $|\Omega|\gtrsim r^2\log^2n$ \cite{CWW:ACHA:18}.

\subsubsection{Robust principal component analysis}
Assume we are given a sum of of a low rank matrix $\BX$ and a sparse matrix $\bm{Y}$:
\begin{align*}
\BD = \BX + \bm{Y}.
\end{align*}
The goal in robust principal component analysis (RPCA) is to reconstruct $\bm{X}$ and $\bm{Y}$ simultaneously from $\bm{D}$. RPCA
appears in a wide range of applications, including video and voice background subtraction  \cite{li2004statistical,huang2012singing}, sparse graphs clustering \cite{chen2012clustering}, 3D reconstruction \cite{mobahi2011holistic}, and fault isolation \cite{tharrault2008fault}. Compared with traditional PCA which computes a low rank approximation to
a data matrix, RPCA is less sensitive to outliers since it includes a sparse component in its formulation. 
RPCA  can be explicitly formulated as 
\begin{equation} \label{eq:rpca}
\min_{\bm{Z},\bm{S}\in\mathbb{R}^{m\times n}} \|\bm{D}-\bm{Z}-\bm{S}\|_F \quad \textnormal{ subject to }\rank(\bm{Z}) \leq r \textnormal{ and } \|\bm{S}\|_0 \leq |\Omega|,
\end{equation}
where $r$ denotes the rank of the underlying low rank matrix $\BX$, $\Omega$ denotes the support set of the underlying sparse matrix $\bm{Y}$, and $\|\bm{S}\|_0$ counts the number 
of non-zero entries in $\bm{S}$.

In \cite{NNSAJ:NIPS:14}, a non-convex algorithm of alternating projections, namely
AltProj,  has been proposed  for \eqref{eq:rpca},
\begin{align*}
\begin{cases}
\BZ_{k+1} = \T_r(\BD-\BS_k)\\
\BSS_{k+1} = \H_{\zeta_{k+1}}(\BD-\BZ_{k+1}).
\end{cases}\numberthis\label{eq:altproj}
\end{align*}
In each iteration, AltProj first computes a new estimate $\BL_{k+1}$ of the low rank component by projecting $\BD-\BSS_{k}$ onto the rank-$r$ matrix manifold $\M_r$ via $\T_r$, and then computes a  new estimate $\BSS_{k+1}$ of the sparse component by 
 projecting $\BD-\BZ_{k+1}$ onto the space of sparse matrices via the entrywise thresholding operator $ \H_{\zeta_{k+1}}$ which is defined by 
 \begin{align*}
 [\H_{\zeta_{k+1}}(\BZ)]_{ij} = 
 \begin{cases}
 Z_{ij} & \mbox{if } |Z_{ij}|>\zeta_{k+1}\\
 0 & \mbox{otherwise}.
 \end{cases}
 \end{align*}
 Here the thresholding value $\zeta_{k+1}$ is adjusted adaptively in each iteration \cite{NNSAJ:NIPS:14}. 
 
Noticing that in the first step of AltProj the SVD on an $n\times n$ matrix is needed to compute the best low rank approximation, we can apply the same idea as in RGrad to reduce the computational cost. That is, before truncating $\BD-\BSS_k$ to its best rank-$r$ approximation, we can first project it onto the tangent space of $\M_r$ at the previous low rank estimate, which leads to the algorithm of accelerated alternating projections (AccAltProj) in \cite{CaiCaiWei}:
\begin{align*}
\begin{cases}
\BZ_{k+1} = \T_r\P_{T_k}(\BD-\BSS_k)\\
\BSS_{k+1} = \H_{\zeta_{k+1}}(\BD-\BZ_{k+1}),
\end{cases}\numberthis\label{eq:accaltproj}
\end{align*}
where $T_k$ is the tangent space of $\M_r$ at $\BZ_k$. Notice that the thresholding values for $\zeta_{k+1}$ in \eqref{eq:altproj} and \eqref{eq:accaltproj} are usually different with each other \cite{NNSAJ:NIPS:14,CaiCaiWei}.

As a result of the additional tangent space projection, AccAltProj is substantially faster than AltProj. Interested readers are referred to \cite{CaiCaiWei} for empirical comparisons of these two algorithms. Moreover, it is established in \cite{CaiCaiWei} that a variant of AccAltProj with a proper initialization is able to successfully separate the underlying low rank and sparse components provided the number of  nonzero entries of the sparse component is not too large.

\paragraph{Remark} Nuclear norm minimization in Section~\ref{sec:convex} and projected gradient descent based on matrix factorization in Section~\ref{sec:pgd} can also be used for spectrally sparse signal reconstruction and robust principal component analysis. We will not present the details here, but refer the reader to \cite{CLMW:JACM:11,CC:TIT:14,CWW:ArXiv:17,yi2016fast} for comprehensive discussion.

%% file: conclusion.tex
\section{Conclusion and discussion}\label{sec:conclusion}
Low rank model plays an important role for exploiting low dimensional structure in high dimensional problems. In this paper, we provide a partial review on effective and efficient approaches for low rank matrix recovery, including nuclear norm minimization, projected gradient descent based on matrix factorization, and  Riemannian optimization based on the embedded manifold of low rank matrices. Theoretical recovery guarantees  have been provided for these approaches. In order to avoid technical details, theoretical results have been presented in an informal way and interested readers could consult related references for comprehensive discussion. 

We make no attempt to cover every aspect of low rank matrix recovery or conduct extensive numerical experiments to evaluate the empirical performance of various algorithms. In this survey, we mainly focus on three measurement models in low rank matrix recovery: matrix sensing, matrix completion and phase retrieval. There are many other low rank reconstruction problems that are not covered, for example low rank matrix demixing \cite{StrohmerWei17}, blind deconvolution \cite{ARR:TIT:14,LLSW:ArXiv:16}, blind demixing \cite{LingStro17}, rank-$1$ measurement model for general low rank matrices \cite{CCGold15,KRTersti14}, and one bit matrix completion \cite{onebitMC}. Recovery guarantees of the algorithms  have been presented  for the noiseless setting. For statistical perspectives in the noisy case, we refer the reader to \cite{CP:XXX:10,CW:ArXiv:15,NegWain12,KolLouTsy11} and references therein for details.